\documentclass[12pt]{amsproc-MAA}
\usepackage{amsfonts,graphicx}

\def\ThmNum{yes}

\def\N{\mathbb {N}}

\def\Z{\mathbb {Z}}
\def\R{\mathbb {R}}
\def\C{\mathbb {C}}

\def\Cstar{\C^*}
\def\disk{\mathbb {D}}

\def\half{\mathbb {H}}

\newfont{\script}{eusm10 at 12pt}
\newfont{\scriptsmall}{eusm10 at 9pt}

\def\ovl{\overline}
\def\phi{\varphi}
\def\eps{\varepsilon}
\def\theta{\vartheta}
\def\Re{\mbox{\rm Re}}
\def\Im{\mbox{\rm Im}}
\def\id{\mbox{\rm id}}

\def\sm{\setminus}

\ifx\ThmNum\undefined 
   \newtheorem{theorem}{Theorem}[section]
   
\else 
   \newtheorem{theorem}{Theorem} 
   
\fi

\newtheorem{lemma}{Lemma}

\newtheorem{corollary}{Corollary}

\def\proof{\par\medskip\noindent {\it Proof. }}
\def\proofof #1 {\par\medskip\noindent {\it Proof of #1. }}
\def\sketch{\par\medskip\noindent{\it Sketch of proof. }}
\def\sketchof #1 {\par\medskip\noindent {\sc Sketch of proof of #1. }}
\def\Box{\framebox[10pt]{\rule{0pt}{3pt}}}
\def\nix{\rule{0pt}{2pt}}
\def\qed{\qedd\par\medskip\noindent}
\def\qedd{\nix\nolinebreak\hfill\hfill\nolinebreak$\Box$}
\def\remark{\par\medskip \noindent {\bf Remark. }}
\def\lineclear
    {\rule{0pt}{0pt}\nopagebreak\par\nopagebreak\noindent}

\newlength{\captionwidth}
\def\reminder #1 {{\sf #1}}
\def\hide #1 {}
\long\def\longhide #1 {}

\hide{
\renewcommand{\section}[1]
{\addtocounter{section}{1}
\medskip \noindent
{\sc \arabic{section}. #1}
\let\@currentlabel\thesection
}
}
\hide{
}

\renewcommand{\lineclear}{\quad}

\newcommand{\newbibitem}[1]
{\bibitem}        % use this version for references like [11]

\def\s{{\underline s}}

\def\Sym{\mathcal{S}}
\def\diam{{\mbox{\rm diam}}}

\begin{document}

\title[Hausdorff Dimension]{Hausdorff Dimension, Its Properties,
and Its Surprises}
\author{Dierk Schleicher}
\address
{Jacobs University Bremen (formerly: International University Bremen),
Research I, Postfach 750 561, D-28725 Bremen, Germany,
{\tt dierk$@$jacobs-university.de}}
%\email{dierk@iu-bremen.de}


\date{\today}

\hide{
\begin{abstract}
We review the motivation and fundamental properties of
the Hausdorff dimension of metric spaces and illustrate this with
a number of examples, some of which are expected and well-known.
We also give examples where the Hausdorff dimension has some
surprising properties: we construct a set $E\subset\C$ of
positive planar measure and with dimension $2$ such that each
point in $E$ can be joined to $\infty$ by one or several curves
in $\C$ such that all curves are disjoint from each other and
from $E$, and so that their union has Hausdorff dimension $1$. We
can even arrange things so that each point in $\C$ which is not
on one of these curves is in $E$. These examples have been
discovered very recently; they arise quite naturally in the
context of complex dynamics, more precisely in the iteration
theory of simple maps such as $z\mapsto \pi\sin(z)$.
\end{abstract}
}

\hide{
\renewcommand{\thefootnote}{}
\footnotetext{{\em Date:} \today}
\renewcommand{\thefootnote}{\arabic{footnote}}
}

\maketitle

%\tableofcontents

%% Attention: fiddling with formatting!
\addtolength{\leftmargini}{-8mm}

\MAAsection{INTRODUCTION}
\label{SecIntro}
The concept of {\em dimension} has many aspects and meanings
within mathematics, and there are a number of very different
definitions of what the dimension of a set should be. The
simplest case is that of $\R^d$: in order to distinguish points
in $\R^d$, we need $d$ different (real) coordinates, so $\R^d$
has dimension $d$ as a (real) vector space. Similarly, a
$d$-dimensional manifold is a space that locally looks like a
piece of $\R^d$.

Another interesting concept is the {\em topological dimension}
of a topological space: every discrete set has topological dimension $0$ (e.g., any finite sets of points in $\R^d$), an injective curve has topological dimension $1$, a disk has dimension $2$ and so on.
The idea is that a set of dimension $d$ can be disconnected in a neighborhood of {\em every} point
by a set of dimension $d-1$: curves and circles can be disconnected by removing isolated points, disks can be disconnected by removing curves and circles, etc.
\hide{\ (the ``at every point'' condition is needed so that the dimension of a set cannot be decreased by enlarging the set: attaching a $1$-dimensional curve
to a $2$-dimensional disk makes it easier to disconnect the set at {\em some} points,
but the topological dimension should still remain $2$).
}
 A formal definition is recursive, starting conveniently with the empty set:  
{\em the empty set has topological dimension $-1$, and a set has topological dimension
at most $d$ if each point has a basis of open neighborhoods whose 
boundaries have topological dimension at most $d-1$.}

\hide{
Therefore, every discrete set has
topological dimension $0$ (e.g., any finite sets of points in $\R^d$):
it can be disconnected at every point by the empty set.
Every simple curve has topological dimension $1$ because it
can be disconnected at every point by removing finitely many points. A disk in the plane
has topological dimension $2$ because it cannot be disconnected by discrete points, but it can be disconnected by lines and circles. The ``at every point'' condition is needed so that the dimension
of a set cannot be decreased by enlarging the set: attaching a $1$-dimensional curve
to a $2$-dimensional disk makes it easier to disconnect the set at {\em some} points,
but the topological dimension should still remain $2$. 
}
\hide{One of several formal definitions of
topological dimension thus goes as follows: }
\hide{; a set has topological dimension exactly $d$
if it has dimension at most $d$ but not at most $d-1$
}

\thispagestyle{plain}

\begin{figure}[p]
\includegraphics[width=50mm]{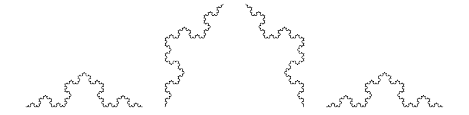}
(b)\!\!\!\!\!\!\!\!\!\!\!\!\!\!\!
\includegraphics[width=50mm,viewport=060 -200 2700 1000,clip]{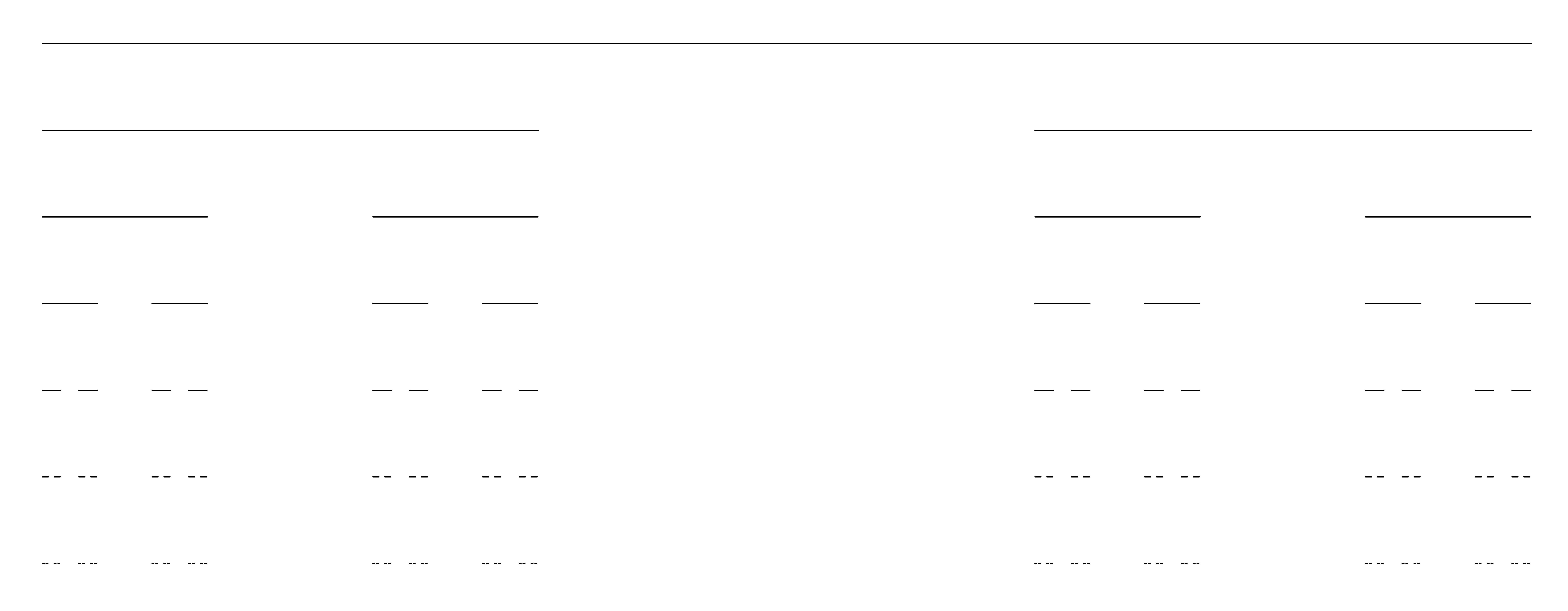} 
\hskip 5mm

(a)\!\!\!\!\!\!\!\!\!\!\!\!\!\!\!
\includegraphics[width=50mm]{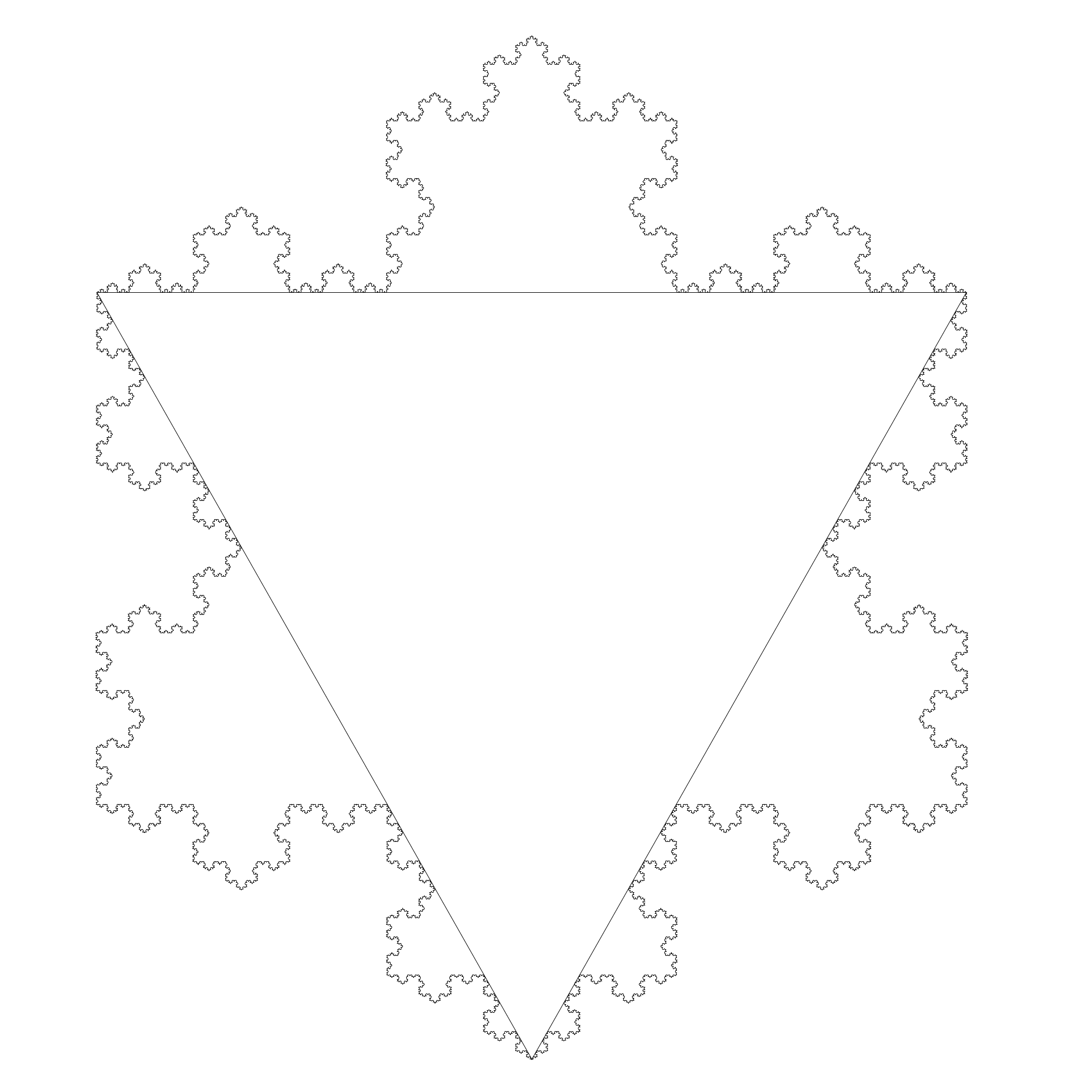}
(c)\!\!\!\!\!
\includegraphics[width=50mm]{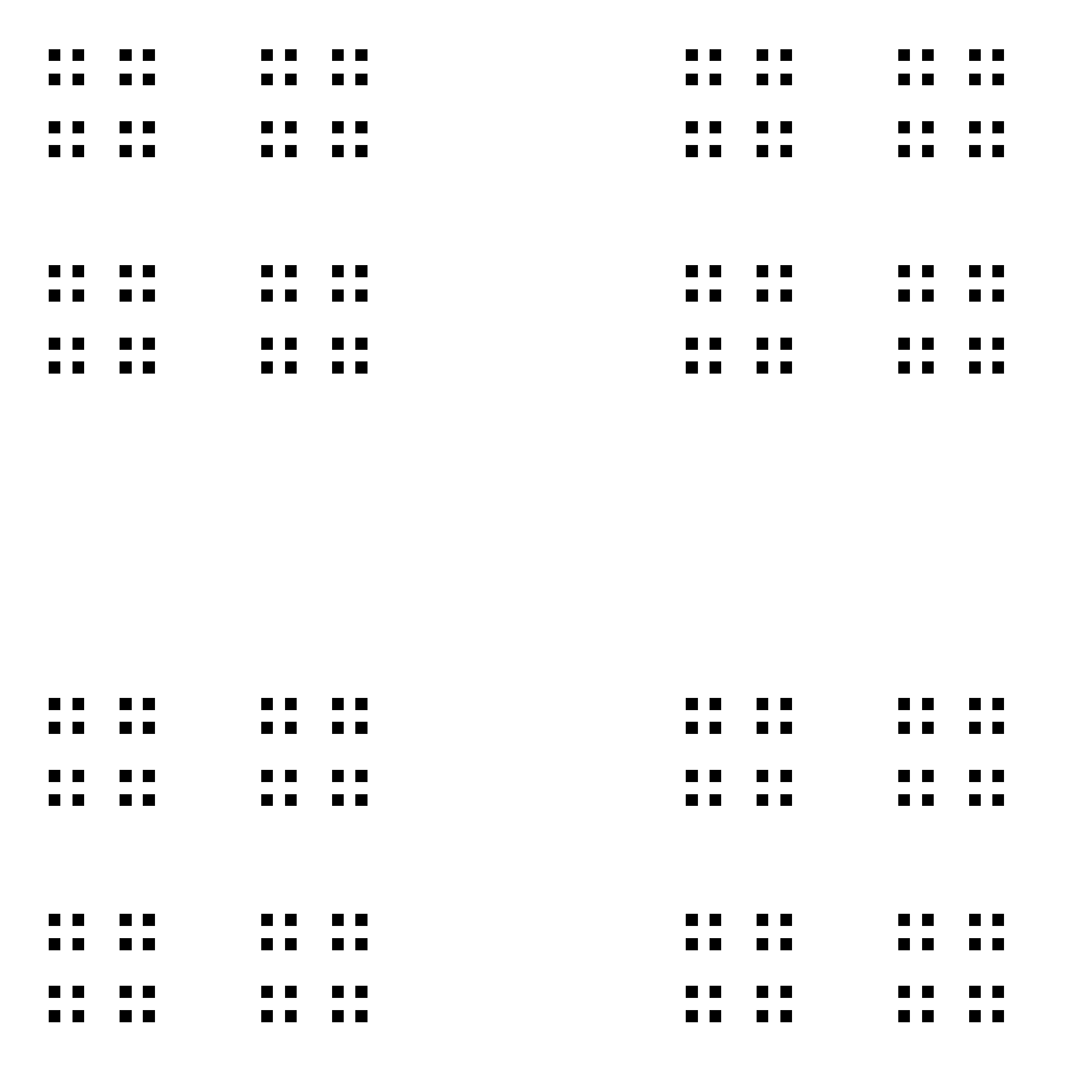} 
\hskip 5mm
\includegraphics[width=50mm]{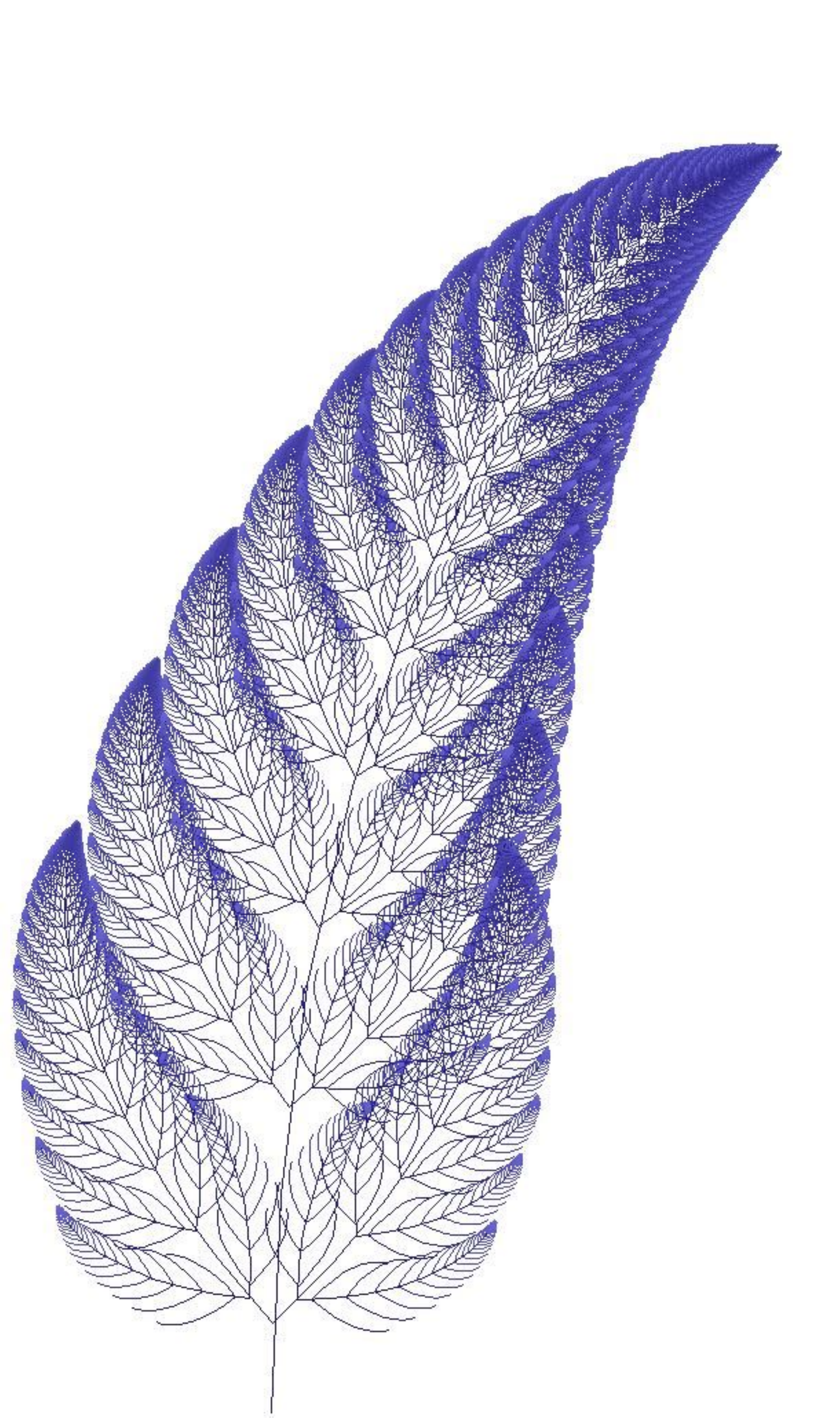}
\includegraphics[width=50mm]{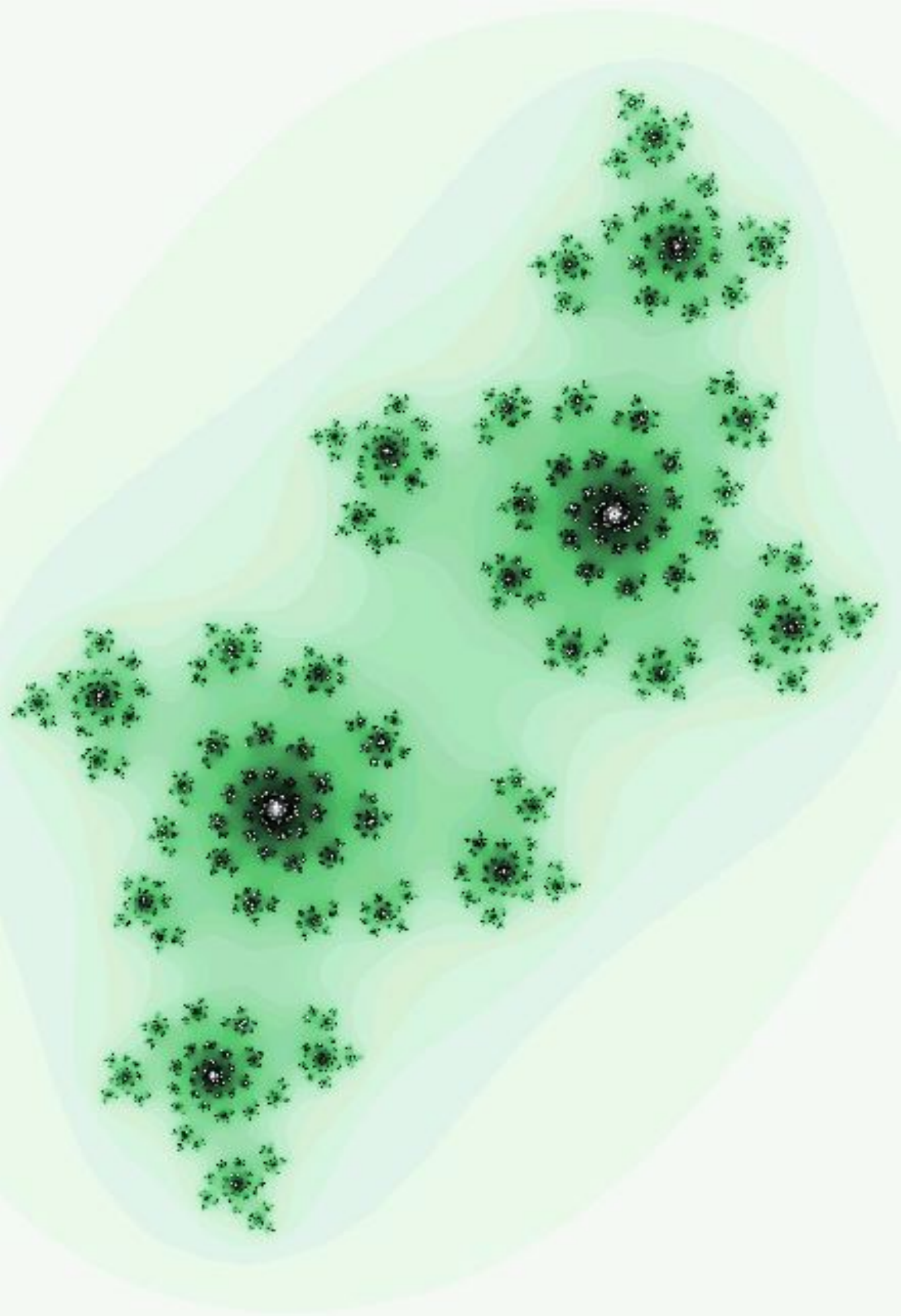}
\par
(d)\qquad\qquad\qquad\qquad\qquad(e)
%\includegraphics[width=50mm]{JuliaCantor.ps}
%\vspace{150mm}\vfill\vfill
\caption{Several ``fractal'' subsets of $\R^2$: 
%{\sl upper left:} 
(a) the ``snowflake'' (von Koch) curve: each of its
three (fractal) sides can be disassembled into four pieces, each of
which is a copy of the entire side, shrunk by a factor $1/3$;
%{\sl upper right:} 
(b) the Cantor middle-third set, consisting of two
copies of itself, shrunk by $1/3$;
%{\sl middle right:} 
(c) a ``fractal'' square in the plane, consisting
of four shrunk copies of itself with a factor $1/3$ (so it has the
same dimension as the snowflake!);
%{\sl lower left:} 
(d) a fern;
%{\sl lower right:} 
(e) the Julia set of a quadratic polynomial. 
}
\label{FigFractals}
\end{figure}

All these dimensions, if finite, are integers (we will ignore
infinite-dimensional spaces). An interesting discussion of
various concepts of dimension, different in spirit from ours, can
be found in the recent article of Manin~\cite{Manin}. 

\looseness +1
We will be concerned with a different aspect of dimension, having
to do with self-similarity of ``fractal'' sets such as those shown in
Figure~\ref{FigFractals}. As 
Mandelbrot points out \cite[p.~1]{Mandelbrot}, ``clouds are not
spheres, mountains are not cones, coastlines are not circles, and
bark is not smooth, nor does lightning travel in a straight
line,'' so many objects occurring in nature are not manifolds.
For instance, 
the fern in Figure~\ref{FigFractals} is constructed by a simple
affine self-similarity process, and people have tried to describe
the hairy systems of roots of trees or plants in terms of
``fractals'', rather than as smooth manifolds. Similar remarks apply
to the human lung or to the borders of most states and countries.

The concept of Hausdorff dimension is almost a century old, but
it has received particularly prominent attention since the advent
of computer graphics and the computer power to simulate and
visualize beautiful objects with importance in a number of
sciences. Earlier, such sets were often constructed by ad hoc
methods as counterexamples to intuitive conjectures. In
the first part of this paper, we try to convince interested
readers that Hausdorff dimension is the ``right'' concept to
describe interesting properties of a metric set $X$: for each
number $d$ in $\R_0^+$ we define the
$d$-dimensional Hausdorff measure $\mu_d(X)$; if $d$ is a positive integer
and $X=\R^d$,
then this measure coincides with Lebesgue measure (up to a
normalization factor). There is a threshold value for $d$,
called $\dim_H(X)$, such that $\mu_d(X)=0$ if $d>\dim_H(X)$ and
$\mu_d(X)=\infty$ if $d<\dim_H(X)$. This value $\dim_H(X)$ is
the {\em Hausdorff dimension} of $X$.

\pagestyle{plain}

We first help to develop intuition for this natural
concept, and then we challenge it by describing a number of
relatively newly discovered sets with very remarkable and
surprising (possibly counterintuitive!) properties of Hausdorff
dimension.
To describe such sets, imagine a curve
$\gamma\colon(0,\infty)\to\C$ that connects the point $0$ to
$\infty$ (we identify a curve $\gamma\colon I\to\C$ with its
image set $\{\gamma(t)\colon t\in I\}$ in $\C$). 
Curves have dimension at least $1$, possibly more, but
the two endpoints certainly have dimension $0$. Now take a
collection of disjoint curves $\gamma_h$, each connecting a
different point $z_h$ to $\infty$. For example, let $z_h=ih$ for $h$ in $[0,1]$ and
$\gamma_h(t)=ih+\gamma(t)$ (provided $\gamma$ is such that all
$\gamma_h$ are disjoint). Then the endpoints are an interval with
dimension $1$, while the union of all curves $\gamma_h$ covers an
open set of $\C$ and should certainly have dimension $2$. This is
true and intuitive: the union of the endpoints has smaller
dimension  than the union of the curves. In this paper, we
describe the following situation \cite{Cosinus}:

\begin{theorem}[A Hausdorff Dimension Paradox]
\label{ThmParadox} \lineclear
There are subsets $E$ and $R$ of $\C$ with the following properties:
\begin{enumerate}
\item
$E$ and $R$ are disjoint;
\item
each path component of $R$ is an injective curve {\em(}a ``ray''{\em)}
$\gamma\colon(0,\infty)\to\C$ connecting some point $e$ of $E$ to
$\infty$ {\em(}i.e., $\lim_{t\to 0}\gamma(t)=e$ and $\lim_{t\to
\infty}\gamma(t)=\infty${\em);}
\item
each point $e$ of $E$ is the endpoint of one or several curves in $R$;
\item
the set $R=\bigcup \gamma((0,1))$ of rays has Hausdorff
dimension $1$;
\item
the set $E$ of endpoints has Hausdorff dimension $2$ and even
full $2$-dimensional Lebesgue measure {\em (}i.e., the set $\C\sm E$ has measure zero{\em);}
\item
stronger yet, we have $E\cup R=\C$: the set of endpoints $E$
is the complement of the $1$-dimensional set $R$, yet each
point in $E$ is connected to $\infty$ by one or several curves in
$R$!
\end{enumerate}
\end{theorem}
\pagestyle{plain}

Mathematics is full of surprising phenomena, and often 
very artful methods are used to construct sets that exhibit these phenomena.
This result is another illustration that many of these phenomena arise quite
naturally in dynamical systems, especially complex dynamics. 
It comes at the end of a series of successively
stronger results. The story started with a surprising result by
Karpi\'nska~\cite{Bogusia}: she established the existence of natural
sets $E$ and $R$ arising in the dynamics of complex
exponential maps $z\mapsto\lambda e^z$ for certain values of 
$\lambda$, where $E$ and $R$ enjoy
properties (1)--(4), as well as (5) in the form that $E$ has Hausdorff
dimension $2$. In \cite{Escaping}, this result was
extended to exponential maps with $\lambda$ in $\Cstar=\C\sm\{0\}$ arbitrary. In
\cite{Guenter}, this was carried over to maps of the form
$z\mapsto ae^z+be^{-z}$; in this case, $E$ always has positive
$2$-dimensional Lebesgue measure.
Finally, condition (6) was established for maps like
$z\mapsto\pi\sinh z$ \cite{Cosinus}.

We start this paper with a
discussion of several concepts of dimension (section~\ref{SecDimension}). 
In section~\ref{SecHausdorff}, we give the definition of Hausdorff
dimension together with a number of its fundamental properties. In
section~\ref{SecKarpinska}, we describe a beautiful example
constructed by Bogus{\l}awa Karpi\'nska in which $E$ has positive
$2$-dimensional Lebesgue measure. In the remainder of the paper, we
show that sets $E$ and $R$ satisfying all the assertions of
Theorem~\ref{ThmParadox}, including $\C=E\dot\cup R$, appear
naturally in complex dynamics, when iterating maps such as
$z\mapsto \pi\sin z$.

The basic features of iterated complex $\sin$ and $\sinh$ maps are
described in section~\ref{SecSinhDynamics}, and a fundamental lemma
for estimating Hausdorff dimension is given in
section~\ref{SecParabola}. In section~\ref{SecFineStructure}, we
then describe the dynamics of the map $z\mapsto \pi\sin z$ in
detail and finish the proof of Theorem~\ref{ThmParadox}. 
Finally, we discuss some related known results about planar Lebesgue
measure, including a theorem of McMullen and a conjecture of
Milnor.

The purpose of this paper is to highlight interesting phenomena
that are observed at the interface between dimension theory and transcendental
dynamics. It cannot serve as an exhaustive survey on the exciting work that
has been done on these two areas, and we can mention only a few
of the most interesting references. A good survey of
transcendental dynamics is found in
Bergweiler~\cite{WalterSurvey}; some more surprising properties
of exponential dynamics are described in
Devaney~\cite{DevaneyTopology}. The topic of ``curves of escaping
points in transcendental dynamics'' was first raised in 1926 by
Fatou~\cite{Fatou} and taken up more systematically by
Eremenko~\cite{Eremenko}. In the special case of exponential
dynamics, it was first investigated by Devaney and
coauthors~\cite{DGH}, \cite{DK} and completed in
\cite{Escaping}, \cite{FRS}. In more general settings, there are existence
results in \cite{DT}, and the current state of the art can be
found in the recent thesis of Rottenfu{\ss}er
\cite{GuenterThesis}. Among current work on Hausdorff dimension in
transcendental dynamics, we would like to mention the survey
papers by Stallard~\cite{Stallard} and by Kotus and
Urba\'nski~\cite{KU}. We apologize to those whose work we have
not mentioned here.

\MAAsection{CONCEPTS OF ``FRACTAL'' DIMENSION}
\label{SecDimension}
The fundamental idea that leads to ``fractal'' dimensions is to investigate interesting sets at different scales of size. Consider a regular three-dimensional cube, say of side-length $1$. We can subdivide this cube into many small cubes of side-length $s=1/k$ for any positive integer $k$. Obviously, the number of little cubes we obtain is $N(s)=k^3=s^{-3}$. However, if we subdivide a unit square into small squares of side-length $1/k$, we obtain $N(s)=s^{-2}$ little squares. The exponent here is the dimension: {\em if a set $X$ in $\R^n$ can be subdivided into some finite number $N(s)$ of subsets, all congruent (by translations or rotations) to one another and each a rescaled copy of $X$ by a linear factor $s$, then the ``self-similarity dimension'' of $X$ is the unique value $d$ that satisfies $N(s)=s^{-d}$, i.e., }
\[
d=\log (N(s))/\log(1/s)
\,\,.
\]

This simple idea can be applied to a number of interesting sets. Consider, for example, the ``snowflake'' curve of Figure~\ref{FigFractals}a: we only look at the top third of the snowflake, above the triangle that we have inscribed for easier description. The detail above the snowflake shows that this top third can be disassembled into $N=4$ pieces, each of which is a rescaled version of the entire top third with a rescaling factor $s=1/3$. The associated dimension must satisfy $3^d=4$ (i.e., $d=\log 4 /\log 3 \approx 1.26\dots$). The snowflake is a curve (and thus has topological dimension $1$), but its self-similarity dimension is greater than that of a straight line: when subdividing a straight line into pieces of one-third the original size, we obtain three pieces; for the snowflake, we get four (and for a square we get nine). Continued refinement has the same dimension: we can break up the four pieces into four pieces each, so that all are rescaled by a factor $s=1/9$; and again $d=\log(4^2)/\log(3^2)=1.26\dots$. 

Let us explore this idea for the standard middle-third Cantor set as shown in Figure~\ref{FigFractals}b. It is constructed by starting with a unit interval, removing the (open) middle third, so as to yield two closed intervals of length $1/3$ each; removing the middle third from these and continuing inductively yields the standard middle-third Cantor set. This set consists of $N=2$ parts (left and right) that both are rescaled versions of the original set with a factor $s=1/3$. This Cantor set has dimension $\log 2/\log 3\approx 0.83\dots$: less than a curve, but more than a discrete set of  points. 

Here is one last example, depicted in Figure~\ref{FigFractals}c: a unit square is subdivided into nine equal subsquares of size $s=1/3$, and only the $N=4$ subsquares at the vertices are kept and further subdivided. The dimension is $\log 4/\log 3\approx 1.26\dots$ as for the snowflake curve. This set is simply the Cartesian product of the middle-third Cantor set with itself.

We can play with the dimension of the Cantor set. For instance, we can start with a unit interval and remove a shorter or longer interval in the middle so as to leave $N=2$ intervals of arbitrary length $s$ in $(0,1/2)$. In the next generations, we always remove an interval in the middle with the same fraction of length, so that the resulting Cantor set is self-similar again. Its dimension is $d=\log 2/\log(1/s)$, and it can assume any real value in $(0,1)$. 

What we have exploited so far is {\em linear self-similarity} of our sets: they consist of a finite number of pieces, each a linearly rescaled version of the entire set. It is only for such sets that the self-similarity dimension applies. Later, we define two further concepts of  ``fractal'' dimension, box-counting dimension and Hausdorff dimension, which make sense for more general sets than the self-similarity dimension; but for the examples we have considered so far, all three dimensions apply and have the same value.

Here is a variation of the construction that leaves the realm of linearly self-similar sets: take the unit interval, replace it with two subintervals of length $s_1\in (0,1/2)$; each of these two intervals is replaced with two further subintervals of length $s_1s_2$ (with $s_2$ in $(0,1/2)$), and so on. If all scaling factors $s_i$ are the same, we have a self-similar Cantor set of dimension $d=\log 2/\log(1/s_i)$ as earlier. If the first $k$ scaling factors are arbitrary, but $s_{k+1}=s_{k+2}=\dots = s$, then our Cantor set consists of $2^k$ small Cantor sets, and these small Cantor sets are linearly self-similar and have dimension $\log 2/\log(1/s)$. If the sequence $s_i$ is not eventually constant, we need a more general concept of dimension. We would expect that the dimension would be $0$ if $s_i\to 0$ and $1$ if $s_i\to 1/2$. This will be true for the box-counting dimension that we define at the end of this section. 

We can even construct a Cantor set within $[0,1]$ that has positive $1$-dimensional Lebesgue measure, so its dimension should certainly be $1$: 
 in the first step, we remove the middle interval of length $1/10$, say; from the
remaining two intervals, we remove the central intervals of length
$1/200$; then we remove four intervals of length $1/4000$, etc..
As a result, the total length of all removed intervals is 
$1/10+2/200+4/4000+\dots=0.1111\ldots=1/9$, so the 
Cantor set left at the end of the process has $1$-dimensional Lebesgue measure $8/9$ (note that
we always remove open intervals, which ensures that the remaining set is
compact, hence has well-defined Lebesgue measure). 

All these Cantor sets are homeomorphic. There is even a
homeomorphism of the unit interval to itself whose restriction to
one Cantor set (say of dimension $0$) yields the other (say of
positive Lebesgue measure). (In general, a nonempty subset of a topological
space is called a {\em Cantor set} if it is compact, totally
disconnected, and without isolated points; any two metric Cantor
sets are homeomorphic \cite[Theorem~2.97]{HY}).

\hide{
To motivate that the dimension of a set can be a non-integer in a natural way, let us consider $d$-dimensional measures $\mu_d$ for all non-negative real numbers $d$. If $d$ is an integer, then $\mu_d$ should coincide with Lebesgue measure in dimension $d$ (perhaps up to a normalization factor): $\mu_1$ should measure length, $\mu_2$ should measure area, and so on. Our measures $\mu_d$ should be outer measures, assigning a value to any bounded set in $\R^n$ (values $0$ and $\infty$ are allowed). For example, if $D_R$ is the disk $\{(x,y)\in\R^2\colon x^2+y^2\le R^2\}$, we want $\mu_2(D_R)=\pi R^2$ and $\mu_1(D_R)=\infty$, while if $C_R=\partial D_R$ is the boundary circle, then we want $\mu_2(C_R)=0$ and $\mu_1(C_R)=2\pi R$. 
Any definition of $\mu_d$ should satisfy the following three properties:
\begin{enumerate}
\item
If $T$ is a rigid motion of $\R^n$ (translation or rotation), then $\mu_d(T(A))=\mu_d(A)$ for every bounded $A\subset\R^n$;
\item
if $S$ is a scaling by a positive real factor $\lambda$, then $\mu_d(S(A))=\lambda^d\mu_d(A)$;
\item
if $A_1,A_2,\dots,A_k$ are disjoint, then $\mu_d(A_1\cup A_2\cup\dots\cup A_k)=\mu_d(A_1)+\mu_d(A_2)+\dots+\mu_d(A_k)$.
\end{enumerate}
In elementary measure theory, one learns that these three properties are (more than) enough to characterize measure of open sets in $\R^n$ (or even Borel sets) uniquely up to a normalization factor: this is the $n$-dimensional measure. However, we want to do much better than that and allow $d$-dimensional measures for all real numbers $d$ in $[0,n]$.
Let us see what the three properties give us in simple cases. If $A$ is a cube in $\R^n$ of side-length $s$, then $A$ is the disjoint union of $N^n$ small cubes $A_1,\dots,A_{N^n}$ of side-length $1/N$, for every integer $N\ge 1$. The first property implies that all little cubes have equal values $\mu_d(A_i)$ for every fixed value of $d$, the second property shows that $\mu_d(A_i)=N^{-d}\mu_d(A)$ for all $i$, and the third property implies that $N^n\mu_d(A_i)=\mu_d(A)=N^{d}\mu_d(A_i)$. The only value of $d$ for which $\mu_d(A_i)$ can have a finite non-zero value is $d=n$.
To discuss self-similarity dimension, let us agree that
$\R$ has dimension $1$. Start with a unit cube in $\R^d$. Shrink
the cube so that its edges now have length $s$ less than $1$ (the
$1$-dimensional size scales with a factor $s$). How many small
cubes are needed to fill the original cube? Obviously we need
$1/s$ cubes in each direction, for a total of $N=(1/s)^d$ cubes.
This exponent $d$ is the (self-similarity) dimension of a cube in $\R^d$:
$d={\log(N)}/{-\log(s)}$. 
With this simple idea, we can determine the self-similarity
dimension of some ``fractals'' like those in
Figure~\ref{FigFractals}: the snowflake curve is a union of three
(fractal) sides. Each side can be decomposed into 
$N=4$ pieces, each of which are scaled images of the entire side,
shrunk by a factor $s=1/3$; therefore, the dimension is
$d=\frac{\log(4)}{-\log(1/3)}=\log(4)/\log(3)\approx 1.26\dots$.
Similarly, if the Cantor middle-third set is shrunk by a factor
of $1/3$, it fills exactly one of its two halves, so we need
$N=2$ copies of the shrunk set to cover the whole:
$d=\log(2)/\log(3)\approx 0.63\dots$. 
Note that we can play with the dimension of the Cantor set: if we
replace the initial interval by two copies, each of length
$s<1/2$, then the dimension is $-\log(2)/\log(s)$ and can be
anywhere in the open interval $(0,1)$. Even values $0$ and $1$
are possible if the factor $s$ is changed in every step: if we
replace the initial unit interval by two intervals of length
$s_1$, and henceforth use a factor $s_2\neq s_1$, then we have
two Cantor sets of dimension $-\log(2)/\log(s_2)$ 
(but their size is different than if we had used the factor $s_2$
to begin with); the entire Cantor set also has dimension
$-\log(2)/\log(s_2)$. Now use different factors
$s_1>s_2>s_3\dots$ in each step; then the resulting Cantor set
should have dimension $\lim_{k\to\infty}-\log(2)/\log(s_k)$ which
can easily become $0$. The same idea can be used to obtain
dimension $1$. Another way to see that dimension $1$ can be
reached is to make the removed intervals in each step
exponentially small: in the first step, we remove not the middle
third, but the middle interval of length $1/10$, say; from the
remaining two intervals, we remove the central intervals of length
$1/200$; then we remove four intervals of length $1/4000$, etc..
As a result, the total length of all removed intervals adds up to
$1/10+2/200+4/4000+\dots=0.1111\dots=1/9$, so the remaining
Cantor set has $1$-dimensional Lebesgue measure $8/9$ (note that
we always remove open intervals, so the remaining set is
compact, hence has well-defined Lebesgue measure). Therefore, as
a set of positive $1$-dimensional Lebesgue measure, the resulting
Cantor set ``must'' have dimension $1$. Note that all these
Cantor sets are homeomorphic to each other, and there is even a
homeomorphism of the unit interval to itself whose restriction to
one Cantor set (say of dimension $0$) yields the other (say of
positive Lebesgue measure). In general, a subset of a topological
space is called a Cantor set if it is compact, totally
disconnected and without isolated points; any two metric Cantor
sets are homeomorphic \cite[Theorem~2.97]{HY}.
}

By taking Cartesian products of linear Cantor sets, we
obtain Cantor subsets of the unit square. We can manufacture these 
so that they have dimension $0$, 
positive $2$-dimensional Lebesgue measure, or anything in
between.

In order to define the dimensions of more general sets like the fern or the
Julia set in Figure~\ref{FigFractals}, we need a more general
approach than self-similarity dimension. For a bounded subset
$X$ of $\R^n$ the idea is as follows: partition $\R^n$ by a
regular grid of cubes of side-length $s$ and count how many of
them intersect $X$; if this number is $N(s)$, then we define the
``box-counting dimension'' (or ``pixel-counting dimension'') of $X$ to be
$\lim_{s\to 0} \log(N(s))/\log(1/s)$.  For example, if $X$ is a
bounded piece of a $d$-dimensional subspace of $\R^n$, then
$N(s)\approx c (1/s)^d$ and the dimension is $d$. This is what a
computer can do most easily: draw the set $X$ on the screen,
count how many pixels it intersects, then draw $X$ in a finer
resolution and count again\dots. Of course, the limit will not
exist in many cases, 
so the box-counting dimension is not always well-defined.
It is, however, well-defined for the linearly self-similar sets discussed earlier,
and for these the self-similarity dimension and the box-counting dimension coincide.
Another drawback of box-counting dimension is that every countable dense subset
$X$ of $\R^n$ has dimension $n$, although a countable set
should be very ``small.'' More generally, this concept of
dimension does not behave well under countable unions. The
underlying reason is that all the cubes used to cover $X$ were
required to have the same size. Giving up this preconception leads
to the definition of Hausdorff dimension.

\MAAsection{HAUSDORFF DIMENSION}
\label{SecHausdorff}
Let $X$ be a subset of a metric space $M$.
We define the $d$-dimensional Hausdorff measure $\mu_d(X)$ of $X$ for any
$d$ in $\R_0^+=[0,\infty)$ as follows:
%\begin{equation}
\[
\mu_d(X)=
\lim_{\eps\to 0}
\inf_{(U_i)}
\sum_i (\diam(U_i))^d
\,\,,
\eqno (*)
%\label{Eq:HausdorffMeasure}
\]
\def\Eq:HausdorffMeasure{*}%
%\end{equation}
where the infimum is taken over all countable covers $(U_i)$ of
$X$ such that $\diam(U_i)<\eps$ for all $i$. The idea is to cover
$X$ with small sets $U_i$ as efficiently as possible (thus the
infimum) and to estimate the $d$-measure of $X$ as the sum of the
$(\diam(U_i))^d$. Smaller values of $\eps$ restrict the set of
available covers, so the infimum can only increase as $\eps$
decreases. Therefore, the limit always exists in $\R_0^+\cup\{\infty\}$.
The measure $\mu_d$ is an outer measure on $M$ for which all
Borel sets are measurable. (Can the reader figure out the meaning of $\mu_0(X)$?)

If $d$ is a positive integer and $X$ is a subset of $M=\R^d$ with its Euclidean metric, 
then the $d$-dimensional Hausdorff measure and the $d$-dimensional 
Lebesgue measure of $X$ coincide up to a scaling
constant (a ball in $\R^d$ of diameter $s$ has $d$-dimensional
Hausdorff measure $s^d$). Also, countable sets have Hausdorff measure
$0$ for all $d>0$. 
The dependence of the $d$-dimensional measures is governed by the
following rather simple lemma:
\begin{lemma}[Dependence of $d$-Dimensional Measure]
\label{LemDependenceMeasure} \lineclear 
For any $d$ in $\R_0^+$ the following statements hold:
\begin{enumerate}
\item
If $\mu_d(X)<\infty$ and $d'>d$, then $\mu_{d'}(X)=0$.
\item
If $\mu_d(X)>0$ and $d'<d$, then $\mu_{d'}(X)=\infty$.
\item
For each bounded set $X$ in a given metric space there is a unique
value $d=:\dim_H(X)$ in $\R_0^+\cup\{\infty\}$ such that
$\mu_{d'}(X)=0$ if $d'>d$ and $\mu_{d'}(X)=\infty$ if $d'<d$.
\end{enumerate}
\end{lemma}

\noindent
The first two assertions of the lemma follow directly from the definition of 
Hausdorff measure in (\Eq:HausdorffMeasure), 
and together they imply the third assertion.

The value $\dim_H(X)$  in Lemma~\ref{LemDependenceMeasure} is called 
the {\em Hausdorff dimension of $X$}. The Hausdorff measure $\mu_d(X)$ 
with $d=\dim_H(X)$ may be zero, positive, or even infinite.

A few remarks might help to elucidate this concept. First, the
definition yields upper bounds for the dimension more easily than
lower bounds: to establish an upper bound for the dimension, it
suffices to find an appropriate covering for each $\eps$; to give
lower bounds, it is necessary to estimate all possible coverings. For
example, the Hausdorff dimension is clearly bounded above by the
box-counting dimension (if the latter exists), but the freedom
to use coverings of varying sizes sometimes yields much smaller
Hausdorff dimension (as mentioned earlier, any
countable set has Hausdorff dimension zero). 

As an example, let $X$ be a bounded subset of a $d$-dimensional
subspace of $\R^n$; to fix ideas, say $X$ is a $d$-dimensional cube. For positive
$s$ let $N(s)$ be the number of open Euclidean balls in $\R^n$ of diameter $s$ needed to
cover $X$. Then $N(s)\leq c(1/s)^d$ for some constant $c$, hence
$\mu_{d'}(X)\leq c (1/s)^d s^{d'}=c s^{d'-d}$. As $s\to 0$,
the latter bound tends to $0$ if $d'>d$, so  $\mu_{d'}(X)=0$ when $d'>d$ and thus
$\dim_H(X)\leq d$. It is not hard to see that
coverings of varying sizes would not change the dimension, so
indeed $\dim_H(X)=d$. This example also shows why we need to take the limit
$\eps\to 0$: if $d<\dim_H(X)$, then coverings using
large pieces would seem to be more efficient, whereas the limit $\eps\to 0$ 
implies that $\mu_d(X)=\infty$ as it should be.

The equivalence between Lebesgue and Hausdorff measures implies
that any set in $\R^d$ with finite positive $d$-dimensional Lebesgue measure
has Hausdorff dimension $d$. This is another indication that
Hausdorff dimension is the ``right'' concept.

It might be instructive to see that for linearly self-similar sets as discussed in
section~\ref{SecDimension}, the Hausdorff dimension never exceeds the
self-similarity dimension. Indeed, if $X$ is a bounded self-similar set of diameter $R$
with the property that $X$ is the union of $N$ subsets, each similar to $X$ and scaled
by a factor $s<1$, then $X$ can be covered by $N$ balls of diameter $sR$, or by
$N^2$ balls of diameter $s^2R$, and so on. Since $s<1$, the diameters tend to zero
as $k\to\infty$. According to the definition
in (\Eq:HausdorffMeasure), this sequence of
finite covers of $X$ yields an upper bound for $\mu_d(X)$ of 
$\lim_{k\to\infty}N^k (s^kR)^d=\lim_{k\to\infty} (Ns^d)^k R^d$, and this is zero if 
$Ns^d<1$ or $d>\log N/\log(1/s)$. Therefore, $X$ has Hausdorff dimension at most 
$\log N/\log(1/s)$. As described earlier, upper bounds for Hausdorff dimension are
easier to give than lower bounds. After all, $X$ might well be countable and thus
have Hausdorff dimension $0$, even though it is linearly self-similar.

The following result collects useful properties of Hausdorff
dimension that are not hard to derive directly from the definition.

\begin{theorem}[Elementary Properties of Hausdorff Dimension]
\label{ThmElementaryProperties}
\lineclear
Hausdorff dimension has the following properties:
\begin{enumerate}
\item
if $X\subset Y$, then $\dim_H(X)\leq\dim_H(Y)$;
\item
if $X_i$ is a countable collection of sets with $\dim_H(X_i)\leq
d$, then $\dim_H\left(\bigcup_iX_i\right)\leq d$;
\item
if $X$ is countable, then $\dim_H(X)=0$;
\item
if $X\subset\R^d$, then $\dim_H(X)\leq d$;
\item
if $f\colon X\to f(X)$ is a Lipschitz map, then
$\dim_H(f(X))\le \dim_H(X)$;
\item
if $\dim_H(X)=d$ and $\dim_H(Y)=d'$, then $\dim_H(X\times Y)\ge d+d'$;
\item
if $X$ is connected and contains more than one point, then
$\dim_H(X)\geq 1$; more generally, the Hausdorff
dimension of any set is no smaller than its topological dimension;
\item
if a subset $X$ of $\R^n$ has finite positive $d$-dimensional Lebesgue
measure, then $\dim_H(X)= d$.
\end{enumerate}
\end{theorem}

For linearly self-similar sets, the Hausdorff dimension coincides with the self-similarity dimension. Thus  Hausdorff dimension is {\em not} preserved under homeomorphisms, as we observed in the case of linear Cantor sets in section~\ref{SecDimension}. Indeed, topology and Hausdorff dimension (or measure theory in general) sometimes have a tenuous coexistence. 

\looseness -1
Some people like the word ``fractal''. 
One possibility is to define a set $X$ to be a ``fractal'' if its Hausdorff dimension is not an
integer ($X$ has ``fractal dimension''). The problem with this definition is that, for example, in $\R^d$
one can have a Cantor set whose Hausdorff dimension is an arbitrary real number in
$[0,d]$ (recall our examples). A curve in $\R^d$ can have any dimension in $[1,d]$, and so on.
Why should a curve $X$ in $\R^d$ be a ``fractal'' when its dimension is $1.001$ or $1.999$, but not when its dimension is $2$? A better definition is this: $X$ is a ``fractal'' if its Hausdorff dimension strictly exceeds its topological dimension. 
More information on ``fractal sets'' and Hausdorff dimension can be found in \cite{Falconer}.

\MAAsection{KARPI\'NSKA'S EXAMPLE}
\label{SecKarpinska}
Here we give a beautiful and surprising example due to
Karpi\'nska.

%\begin{Example} [Karpi\'nska]
%\label{ExKarpinska}
%\lineclear
\medskip\noindent
{\bf Example (Karpi\'nska).}
There exist sets $E$ and $R$ in the complex plane $\C$ with the following
properties:
\begin{enumerate}
\item
$E$ and $R$ are disjoint;
\item
$E$ is totally disconnected but has finite positive
$2$-dimensional Le\-bes\-gue measure (hence $E$ has topological
dimension $0$ and Hausdorff dimension $2$);
\item
each connected component of $R$ is a curve connecting a single
point of $E$ to $\infty$;
\item
$R$ has Hausdorff dimension $1$.
\end{enumerate}
%\end{Example}
\medskip

Why is this surprising? Each connected component of $R$ is a
single curve connecting one point of $E$ to $\infty$, so each
connected component of $E\cup R$ contains one point of $E$ and a
whole curve in $R$. The set $E\cup R$ is an uncountable union of such
things, a union so large that the union of all these single points
of $E$ acquires positive $2$-dimensional Lebesgue measure, hence
Hausdorff dimension $2$. In the same union, the dimension of $R$ stays $1$, 
so a $1$-dimensional set can be big enough to
connect each point in the $2$-dimensional set $E$ to $\infty$ via
its own curve, all curves and endpoints being disjoint!

Once this phenomenon is discovered (which happened unexpectedly in
complex dynamics \cite{Bogusia}), its proof is surprisingly
simple. For the set $E$ we use a Cantor set made from an initial closed
square, which is replaced with four disjoint closed subsquares, each of which is
in turn replaced with four smaller disjoint subsquares, etc. 
It is quite easy to arrange the sizes of the squares so that the
resulting Cantor set has positive area: one simply has to make
sure that the area lost at each stage is small enough
so that the cumulative area lost is less than, say, half the area of the initial
square. This leaves a Cantor set with positive area (which is
simply a product of two one-dimensional Cantor sets with positive
$1$-dimensional measure).

\begin{figure}[htbp]
\includegraphics[viewport=095 190 480 620,width=\textwidth,clip]{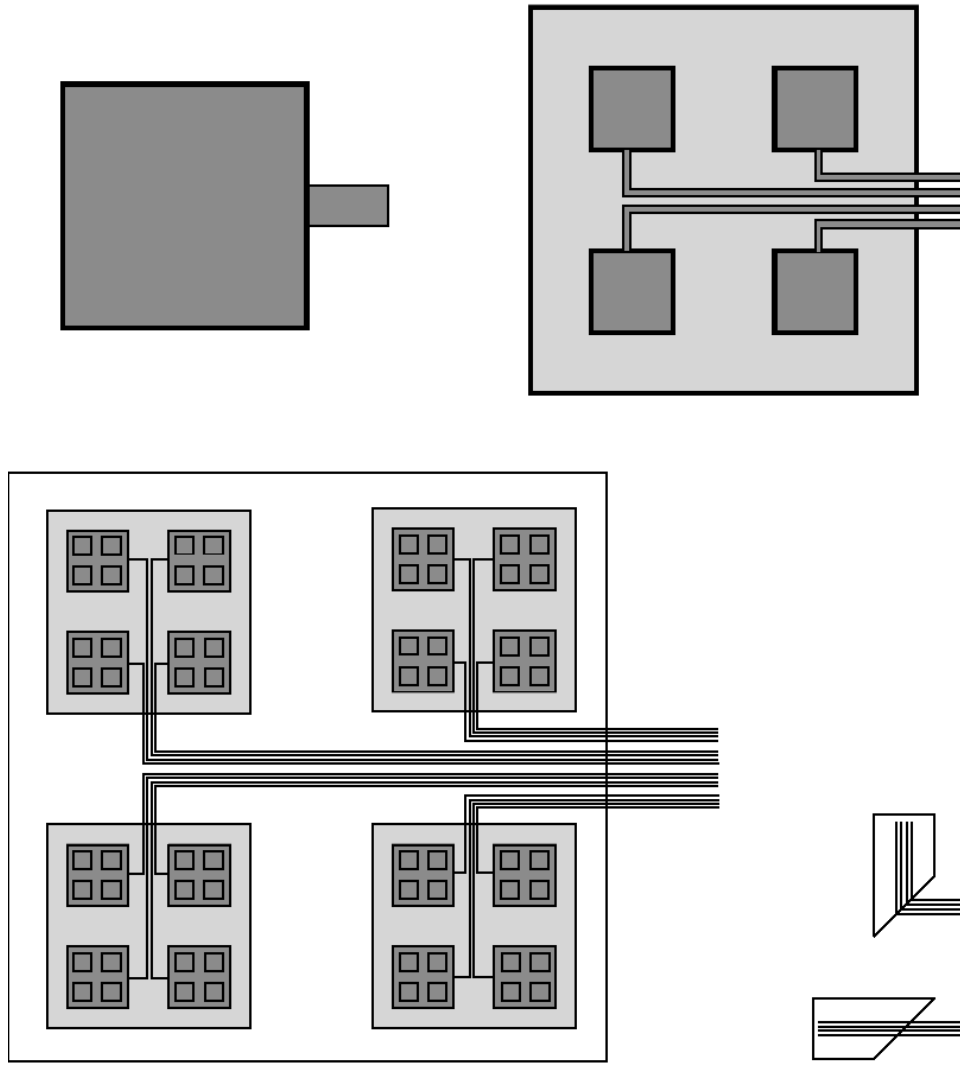}
\caption{The construction of Karpi\'nska's example. Shown are the initial square
and the initial rectangle, as well as two refinement steps. In each step, we keep 
the dark shaded area, so we have a nested sequence of compact sets
(the area of the previous refinement step is shown in a lighter shade). The
detail in the lower right shows that a Cantor set of curves can be
given a right-angled turn by replacing a subset with its
mirror-image, not changing the dimension. 
}
\label{FigKarpinska}
\end{figure}

The construction of the curves is indicated in Figure~\ref{FigKarpinska}.
We start with an initial rectangle that terminates at the initial square. When
the square is refined into four closed subsquares, the rectangle is
subdivided into four parallel closed subrectangles and extended through
the initial square so that the four extended subrectangles reach
the four subsquares. This process can be repeated at each subsequent stage
to create a collection of ``rectangular tubes'' connecting the $4^n$ squares 
in the $n$th subdivision step with the right side of the original square. 
The $n$th refinement step yields $4^n$ squares, each of which has a 
``rectangular tube'' attached to it, so that we have $4^n$ connected components.
Let $X_n$ be the set constructed in step $n$ (consisting of $4^n$ squares together
with their ``rectangular tubes''). Then $X_{n+1}$ is a subset of $X_n$.
More precisely, each step refines each of the $4^n$ connected components of $X_n$ 
into four connected components of $X_{n+1}$. 

It is clear that the countable intersection $\bigcap X_n$
yields a compact set $X$ with the following properties: each connected component
of $X$ consists of one point of $E$ and a curve connecting that point to the
right end of the initial rectangle. Set $R=X\sm E$. All that
remains to show is that $R$ has Hausdorff dimension $1$. Observe
that $R$ restricted to the initial rectangle is a product of an
interval (in the horizontal direction) with a Cantor set (in the
vertical direction). We can arrange things so that the vertical
Cantor set has Hausdorff dimension $0$, so the subset of $R$
within the initial rectangle has Hausdorff dimension $1$. Next
consider the subset of $R$ within the original square but outside
of the first generation subsquares. This looks like a Cantor set
of curves as before, but with a right-angled turn in the middle.
If half of this curve is replaced with its mirror-image,
we obtain a proper Cantor set of curves with dimension $1$ (see
the detail in Figure~\ref{FigKarpinska}), and this reflection does
not change the Hausdorff dimension. The entire set $R$ is a
countable union of such $1$-dimensional Cantor sets of curves,
each with one turn, that become smaller as they approach $E$. Therefore,
$R$ still has dimension $1$.

The last small issue is that the curves in $R$ do not connect $E$
to $\infty$, for they terminate at the right end of the initial
rectangle. This shortcoming can be cured by extending the initial rectangle
to the right by countably many copies of itself.

Certainly, one might find this result surprising. Is it an artifact of the concept of
Hausdorff dimension, indicating that its definition is problematic?
The answer is no: a weaker form of this surprise occurs even from the point
of view of planar Lebesgue measure. Our construction assures that $R$ has
zero planar measure, whereas $E$ has strictly positive planar measure. 
Hausdorff dimension is a way
of making the surprise more precise and stronger; the
surprise lies in the sets $E$ and $R$, not in any definition.

%\looseness -1
We conclude this section with an example of an ``impossible'' set that was brought to our attention by Adam Epstein: Larman~\cite{Larman} defines a compact set in $\R^n$ (for any $n\ge 3$) that is the disjoint union of closed line segments and has positive $n$-dimensional Lebesgue measure. However, removing the two endpoints from each segment, a set with zero measure remains (this is impossible in $\R^2$). In other words, we have a bunch of uncooked spaghetti in $n$-space so that all the nutrition lies in the endpoints. We now proceed to show how much better we can do, even in $\R^2$, when using cooked spaghetti and complex dynamics.

\MAAsection{DYNAMICS OF COMPLEX SINE MAPS}
\label{SecSinhDynamics}
%\looseness -1
In the rest of this article, we describe how a much stronger
result arises quite naturally in the study of very simple dynamical systems,
such as the one given by iterating as simple a map (apparently!)
as $z\mapsto \pi\sinh z$ on $\C$ (but recall that Karpi\'nska
developed her example of section~\ref{SecKarpinska} only after
she had discovered an analogous phenomenon in the dynamics of
exponential maps). We again have sets $E$ and
$R$ as in Karpi\'nska's example, but this time 
$E\cup R=\C$. As before, each path component of $R$ is
a curve connecting one point in $E$ to $\infty$, and $R$ still has
Hausdorff dimension $1$, but now the set $E=\C\sm R$ has infinite Lebesgue
measure, even full measure in $\C$, and is so big that its
complement has dimension~$1$---nevertheless, each point of $E$ can be
connected to $\infty$ by one or even several curves in $R$!

We set up the construction as follows. Let $f\colon\C\to\C$ be
given by $f(z)=k\pi\sinh z=(k\pi/2)(e^z-e^{-z})$ with a nonzero integer
$k$. We study the dynamics given by iteration
of $f$: by $f^{\circ n}$ we denote the $n$th iterate of $f$
(i.e., $f^{\circ 0}=\id$ and $f^{\circ (n+1)}=f\circ f^{\circ n}$).
Of principal interest is the set of ``escaping points,'' meaning the set
\[
I:=\{z\in\C\colon f^{\circ n}(z)\to\infty\mbox{ as
$n\to\infty$}\}
\]
consisting of those points that converge to $\infty$ under iteration of $f$
(in the sense that $|f^{\circ n}(z)|\to\infty$).
Here, $I$ stands for ``infinity''; this set plays a fundamental
role in the iteration theory of polynomials
\cite[sec.~18]{MiIntro} and is just beginning to emerge as
equally important for transcendental entire functions.
\hide{Section~\ref{SecComplexDynamics}. }
Eremenko \cite{Eremenko} has shown that for every transcendental
entire function the set $I$ is nonempty, and he asked whether
every path component of $I$ was unbounded. An affirmative
answer to this question is currently known only for functions of
the form $z\mapsto\lambda e^z$ \cite{Escaping} or $z\mapsto
ae^z+be^{-z}$ \cite{Guenter}, where $\lambda$, $a$, and $b$ are 
nonzero complex numbers. The latter family includes our functions $f$. 
(Recently, this question was answered affirmatively in greater generality
in \cite{GuenterThesis}, \cite{RRRS}, and \cite{Baranski}. However, Eremenko's
question is not true for all transcendental functions; counterexamples are
constructed in \cite{GuenterThesis}, \cite{RRRS}).
The following is a special case of what is known for this family \cite{Cosinus}:
\hide{
\cite{Guenter,Cosinus}%
\footnote{Note that the parametrization of our maps
$g\colon(0,\infty)\to I$ differs from the one used in
\cite{Guenter,Cosinus}} 
}

\begin{theorem}[Dynamic Rays of Sine Functions]
\label{ThmEscapingSine} \lineclear
\begin{enumerate}
\item
\label{ItemThm1}
For the function $f(z)=k\pi\sinh z$ with a nonzero integer $k$ each
path component of $I$ is a curve $g\colon(0,\infty)\to
I$ or $g\colon[0,\infty)\to I$ such that
$\lim_{t\to\infty}\Re\, g(t)=\pm\infty$.
Each curve is contained in a horizontal strip of height $\pi$.
{\em(}These curves are called {``dynamic rays.''}{\em)}
\item
\label{ItemThm2}
For each such curve $g$ the limit
$z:=\lim_{t\searrow 0}g(t)$ exists in $\C$ and is called the {
``landing point''} of $g$ {\em(}``the dynamic ray $g$ lands at $z$''{\em)}. If
$t>t'>0$, then the two points
$g(t)$ and $g(t')$ escape in such a way that 
\[
|\Re f^{\circ k}(g(t))|-|\Re f^{\circ k}(g(t'))|\longrightarrow\infty
\,\,.
\]
\item
\label{ItemThm3}
Conversely, every point $z$ of $\C$ 
either is on a unique dynamic ray or is the landing point of one, two,
or four dynamic rays {\em(}i.e., either $z=g(t)$ for a unique dynamic
ray $g$ and a unique $t>0$, or $z=\lim_{t\searrow 0}g(t)$
for up to four rays $g${\em).} 
\end{enumerate}
\end{theorem}

We will indicate in section~\ref{SecFineStructure} why these
results are not too surprising, even though the precise proofs are
technical. This leads quite naturally to a decomposition
$\C=E\dot\cup R$ as required for our result:
\[
R:=\bigcup_{\mbox{\scriptsize rays $g$}} g((0,\infty))
\,\,,
\qquad
E:=\bigcup_{\mbox{\scriptsize rays $g$}} \lim_{t\searrow 0}
g(t)
\,\,.
\]
If your intuition for the complex sine map is better than for the hyperbolic variant,
then you may use the former instead: the situation is
exactly the same, except that the complex plane is rotated by
$90^0$. We prefer to use the $\sinh$ map because in half-planes far to the left
or far to the right it is essentially the same as $z\mapsto e^{-z}$ and $z\mapsto e^z$,
respectively (up to a factor of $2$). Note also that the parametrization of our rays
$g\colon(0,\infty)\to I$ differs from the one used in
\cite{Guenter} and \cite{Cosinus}.

\MAAsection{THE PARABOLA CONDITION}
\label{SecParabola}
The driving force behind our results is a fundamental lemma of Karpi\'nska~\cite{Bogusia},
adapted to fit our purposes. For real numbers $\xi$ in $(0,\infty)$ and $p$ in $(1,\infty)$ 
consider the sets
\[
P_{p,\xi}:=\left\{x+iy\in\C\colon |x|>\xi, |y|<|x|^{1/p}\right\}
\]
(the ``$p$-parabola,'' restricted to real parts greater than $\xi$).
Also let $I_{p,\xi}$ be the subset of $I$ consisting of those escaping points $z$
for which $f^{\circ n}(z)$ is in $P_{p,\xi}$ for all $n$ (the set of
points that escape within $P_{p,\xi}$). The results in this
section hold for all maps $f(z)=ae^z+be^{-z}$ with
$a$ and $b$ nonzero complex numbers.

\begin{lemma}[Dimension and the Parabola Condition]
\label{LemDimParabola} \lineclear
For each $p$ in $(1,\infty)$ and each sufficiently large $\xi$, the set $I_{p,\xi}$ has Hausdorff dimension at most $1+1/p$.
\end{lemma}
\proof
First observe that we seek only an upper estimate for the
Hausdorff dimension. Therefore it suffices to find a family of covers
whose sets have diameters less than any specified $\eps>0$ so that their
combined $d$-dimensional Hausdorff measure is bounded for each $d$ with
$d>1+1/p$. For bounded subsets of $I_{p,\xi}$, we
construct a finite cover in ``generations'' zero, one, two, {\ldots}
so that each set in the $n$th generation is refined into finitely many
smaller sets in the $(n+1)$th generation. We do this in such a way that the
diameters of all sets tend to zero as the number $n$ of generations tends to infinity, and so
that the combined $d$-dimensional Hausdorff measure of all sets in the
$n$th generation decreases as $n$ tends to infinity provided that $d>1+1/p$.
In view of the definition in (\Eq:HausdorffMeasure), this implies that the 
$d$-dimensional Hausdorff measure of $I_{p,\xi}$ is finite whenever $d>1+1/p$, 
hence that the Hausdorff dimension of $I_{p,\xi}$ is at most $1+1/p$.

We first outline the proof while making a number of
simplifications; we then argue that these do not matter. 
The first simplification is that when $\Re\, z>\xi$, we write
$f(z)=a e^z$ (ignoring the exponentially small error term $b
e^{-z}$), and when $\Re\, z<-\xi$, we write $f(z)=b e^{-z}$. 
For simplicity, we ignore certain bounded factors: we do not
distinguish between side-lengths and diameters of squares, and we
suppress factors like $\pi/|a|$ or $\pi/|b|$ that appear
all over the place but influence only Hausdorff measure,
not dimension.

For the purposes of this proof, ``standard square'' means a closed
square of side-length $\pi$ with sides parallel to the coordinate
axes. The image $f(Q)$ of a standard square $Q$ is a semiannulus bounded by two
semicircles and two straight radial boundary segments. If the imaginary parts
of $Q$ are varied while the real parts are kept fixed, then the semiannulus $f(Q)$
rotates around the origin. We always adjust the imaginary parts of our standard
squares so that $f(Q)$ is entirely contained in the right or the left half-plane, 
which is equivalent to the condition that the two straight radial boundary segments 
of $f(Q)$ are contained in the imaginary axis.

\begin{figure}[htbp]
\includegraphics[viewport=145 450 490 720,width=140mm,clip]{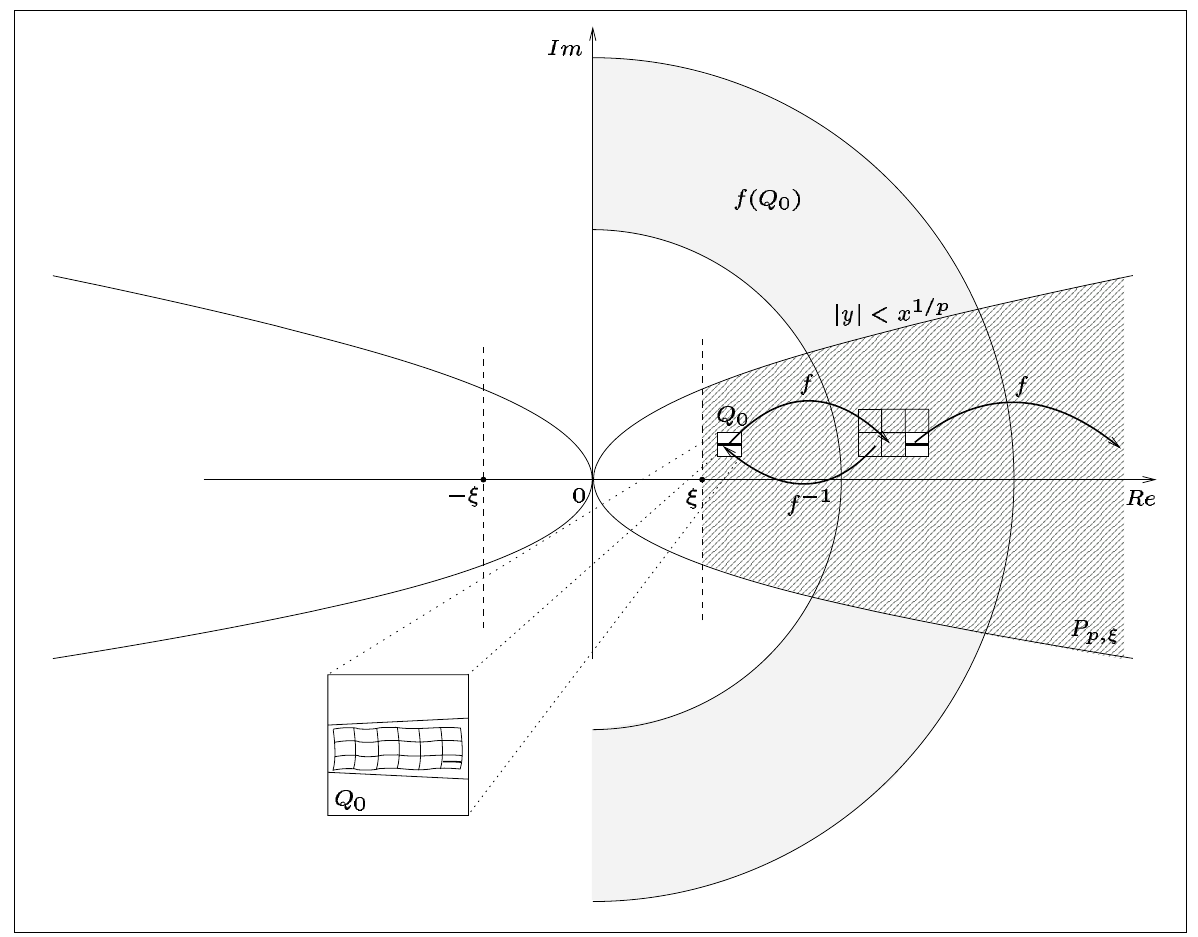}
\caption{Calculating the Hausdorff measure of $I_{p,\xi}$ involves
a partition by iterated preimages of a square grid, as well as
refinements of such a partition.}
\label{FigParabolaDimension}
\end{figure}

Cover $P_{p,\xi}$ by a countable collection of standard squares
with disjoint interiors. 
Fix any particular square $Q_0$ with real parts in
$[x,x+\pi]$, where $x\geq\xi$ and $\xi$ is sufficiently large (the case where
$x\leq-\xi$ is analogous). Now $f(Q_0)$ intersects $P_{p,\xi}$ in
an approximate rectangle with real parts between $\pm|a| e^x$ and
$\pm|a| e^{x+\pi}$ and imaginary parts at most $(|a|
e^{x+\pi})^{1/p}=(|a| e^\pi)^{1/p}e^{x/p}$. Therefore, the number
of standard squares of side-length $\pi$ needed to cover
$f(Q_0)\cap P_{p,\xi}$ is approximately $ce^x\cdot
e^{x/p}=ce^{x(1+1/p)}$, where
\[
c=|a|(e^\pi-1)\cdot 2(|a|e^\pi)^{1/p}/\pi^2=2(e^\pi-1)e^{\pi/p}|a|^{1+1/p}\pi^{-2}
\,\,.
\]
Transporting
these squares back into $Q_0$ via $f^{-1}$, we cover not all of
$Q_0$, but all those points $z$ of $Q_0$ with $f(z)$ in $P_{p,\xi}$
(see Figure~\ref{FigParabolaDimension}). Since $|f'(z)|>|a| e^x$ on
$Q_0$, the covering sets are approximate squares of side-length
at most $(\pi/|a|)e^{-x}$, hence diameter at most $(\sqrt
2\pi/|a|)e^{-x}$. Ignoring bounded factors, we simplify this value to $e^{-x}$. 
We call this covering the ``first generation covering'' within
$Q_0$ (while $\{Q_0\}$ itself is the zeroth generation covering).

Let us see what effect this refinement has on the $d$-dimensional Hausdorff measure.
The covering at generation zero is a standard square and has constant measure.
In generation one, the covering of $Q_0$ has measure
$\sum (\diam(U_i))^d \approx c e^{x(1+1/p)} (e^{-x})^d=ce^{x(1+1/p-d)}$. 
Since $d>1+1/p$, this is small for large $x$ in $(\xi,\infty)$, so this first refinement reduces
the measure.

\hide{
Recall the definition of $d$-dimensional Hausdorff measure of a set $X$ in 
(\ref{Eq:HausdorffMeasure}), using coverings $U_i$ of $X$. In our case,
the first generation covering of $Q_0\cap f^{-1}(P_{p,\xi})$ has
$\sum (\diam(U_i))^d \approx c e^{x(1+1/p)} (e^{-x})^d=ce^{x(1+1/p-d)}$. 
This is not yet the $d$-dimensional Hausdorff measure of $Q_0\cap f^{-1}(P_{p,\xi})$
because we have to consider the infimum over all coverings using sets of diameter
at most $\eps$, and we have to take the limit as $\eps\to 0$. We already
see the crucial term $e^{x(1+1/p-d)}$ which tends to $0$ as $x\to\infty$ when $d>1+1/p$.
}

\hide{
Considering other coverings with diameter less than $\eps$ can only decrease the value
of the sum: in order to show that the dimension is at most $1+1/p$, it suffices to find
a single cover for each $\eps>0$ so that the sum remains finite. We do this by refining our squares. 
}
We continue to refine our coverings so that the diameters of the covering sets tend to zero, while
the $d$-dimensional Hausdorff measure does not increase. 
Each approximate square of generation $n$ gets replaced with some number of
much smaller approximate squares of generation $n+1$. 
What brings the dimension down is that we consider only orbits in
$P_{p,\xi}$, throwing away everything that leaves this parabola
under iteration. We may thus maintain the inductive claim that all
approximate squares of generation $n$ have images under $f$,
$f^{\circ 2}$, \dots, $f^{\circ n}$ that intersect $P_{p,\xi}$; moreover, 
if $Q'$ is an approximate square of generation $n$, then
$f^{\circ n}(Q')$ is a standard square whose points have very large real parts,
say in $[y,y+\pi]$ for some  $y$ satisfying $y\ge \xi$. 

Let $\lambda:=|(f^{\circ n})'(z)|$ for some
$z$ in $Q'$ (this derivative is essentially constant on $Q'$, as noted later). 
Then $Q'$ is an approximate square of side-length
$\pi/\lambda$, so it contributes approximately $\pi^d/\lambda^d$
to the $d$-dimensional Hausdorff measure. We now determine
what happens to this measure under refinement.

Just as in the first step,
$f^{\circ(n+1)}(Q')\cap P_{p,\xi}$ is covered by
$N_y:=ce^{y(1+1/p)}$ standard squares of side-length $\pi$, so
the standard square $f^{\circ n}(Q')\cap f^{-1}(P_{p,\xi})$ is
covered by $N_y$ approximate squares of side-length
$(\pi/|a|)e^{-y}$ or $(\pi/|b|)e^{-y}$. Ignoring constants again, we simplify this
to $e^{-y}$. We need $N_y$ very small approximate
squares to cover those points in $Q'$ that remain in $P_{p,\xi}$
for $n+1$ iteration steps. 
These $N_y$ approximate squares within $Q'$ have side-lengths
approximately $e^{-y}/\lambda$, so their contribution to
the $d$-dimensional Hausdorff measure within $Q'$ is roughly $N_y\cdot
(e^{-y}/\lambda)^d =ce^{y(1+1/p-d)}\lambda^{-d}$, whereas
the contribution of $Q'$ before refinement was $\pi^d\lambda^{-d}$.
Therefore, if $d>1+1/p$, each refinement step reduces the
$d$-dimensional Hausdorff measure (at least when $\xi$ is large).
It follows that the $d$-dimensional Hausdorff measure of $Q_0\cap I_{p,\xi}$ is finite
whenever $d>1+1/p$, so
Lemma~\ref{LemDependenceMeasure} implies that 
\[
\dim_H(Q_0\cap I_{p,\xi})\leq 1+1/p
\,\,.
\]
Since $I_{p,\xi}$ is a countable union of sets of dimension at most $1+1/p$, the claim follows.

There are two main inaccuracies in this proof: we have ignored
constants, and we have ignored the geometric distortions caused by the
mapping $f$ and its iterates. The latter are induced by two
problems: we have disregarded one of the two exponential terms in
$f$, and the continued backward iteration of standard squares
under a finite iterate of $f$ might distort the shape of the
squares because $f'$ or $(f^{\circ n})'$ is not exactly constant
on small approximate squares. However, this distortion problem
is easily cured by a useful lemma usually called the {\em Koebe
Distortion Theorem} \cite[Theorem~2.7]{McBook} for conformal mappings: {\em  for
$r\ge 1$ let $\disk_r:=\{z\in\C\colon |z|<r\}$, and let 
$K_r$ be the family of injective holomorphic mappings
$g\colon\disk_1\to\C$ that have extensions to $\disk_r$
as injective holomorphic mappings. Then for each $r>1$ all maps $g$ in
$K_r$ have distortions (on $\disk_1$) that are uniformly bounded in terms only
of $r$}. Here the precise definition of distortion is
irrelevant: any quantity can be used that measures the deviation
of $g$ from being an affine linear map. A more precise way of
stating this result is as follows: if we normalize so that
$g(0)=0$ and $g'(0)=1$, then the space $K_r$ is compact (in the
topology of uniform convergence). You may want to remember this
fact as the ``yellow of the egg theorem'': when you spill an egg
into a frying pan, the whole egg can assume any shape (this represents the Riemann
map from the disk of radius $r>1$ onto a simply connected domain in
$\C$), but its smaller yolk (the yellow of the egg, represented by the unit disk) is not
distorted too much (it remains essentially a round disk, and derivatives
at any two points differ at most by a bounded factor).

In our context, the maps are easily seen to have bounded
distortion, so we may assume that the $n$th iterate $f^{\circ n}$,
which maps an $n$th generation approximate square to a standard
square, is a linear map with constant complex derivative. All this
does is to introduce a bounded factor in the diameters and in the
number of sets in the coverings. These factors do not increase under 
repeated refinement.

The second simplification was that at several stages we ignored certain
bounded factors. For example, in the calculation
of Hausdorff measures, we replaced diameters with
side-lengths. This introduces a factor of $\sqrt 2$ into the measure
estimates, but it has no impact on the dimension.
Similarly, we have ignored factors like $\pi/|a|$ or $\pi/|b|$,
we have counted the number of necessary squares only
approximately, ignoring boundary effects, and we have
assumed that the derivative of $f^{\circ n}$ is constant on
small approximate squares. Each of these simplifications might lead to a change in
the Hausdorff measure by a bounded factor, but the dimension remains unaffected.
The crucial fact is that refinements do not increase the $d$-dimensional measure
when $d>1+1/p$ and $x$ is sufficiently large, and this fact is correct.
\qed

We have now shown that escaping orbits that spend their entire
lives within the truncated parabolas $P_{p,\xi}$ form a very
small set. It is easy to see that the same is true for the set of points
that spend their entire orbits within $P_{p,\xi}$ except for
finitely many initial steps (see Corollary~\ref{CorHausdorffRays}).
Nonetheless, the surprising fact is that from a different
(topological) point of view, {\em most} orbits do exactly
that: after finitely many initial steps, they enter
$P_{p,\xi}$ and remain there. All this is based on the following
result.

\begin{lemma}[Horizontal Expansion]
\label{LemHorizExpansion} \lineclear
\hide{Fix a map $f(z)=ae^z+be^{-z}$ with $ab\neq 0$.}
For each $h>0$ there is an $\eta>0$ with the following
property: if $(z_k)$ and $(w_k)$ are two orbits such that
$|\Im(z_k-w_k)|<h$ for all $k$ and $|\Re\, z_1|>|\Re\, w_1|+\eta$,
then for each pair $p$ and $\xi$ there is an $N$ such that $z_k$
belongs to $P_{p,\xi}$ whenever $k\geq N$.
\end{lemma}
\sketch
We do not give a precise proof, which involves easy but lengthy
estimates. Instead, we outline the main idea, again ignoring bounded
factors. Let $c:=\max\{|a|,|b|\}$ and $c':=\min\{|a|,|b|\}$, where
$f(z)=ae^z+be^{-z}$. We start by estimating $\Re f(w)$ for sufficiently
large $|\Re\, w|$:
\[
|\Re f(w)|+c\leq |f(w)|+c \leq 
c \exp|\Re\, w|+c < \exp(|\Re\, w|+c)
\,\,,
\]
which yields $|\Re\, w_{k+1} |\le |w_{k+1}|<\exp^{\circ k}(|\Re\, w_1|+c)$
by induction. Therefore
\[
|\Im\, z_{k+1} |\leq
|\Im\, w_{k+1} |+h\leq
|w_{k+1}|+h\leq
\exp^{\circ k}(|\Re\, w_1 |+c)+h
\,\,.
\]

If $|\Re\, z|>|\Re\, w|+\eta$ and both are sufficiently large, then
\[
|f(z)|\geq c'\exp|\Re\, z|>c'\exp(|\Re\, w|)\exp\eta\approx|f(w)| e^\eta
\,\,,
\]
hence $|f(z)|\gg |f(w)|$ if $\eta$ is large. Since the
imaginary parts of $f(z)$ and $f(w)$ are approximately equal, 
the absolute value of $f(z)$ must come mainly from its real part, so 
\[
|\Re\, f(z)|-1\geq \frac 1 e |f(z)|\approx\exp(|\Re\, z|-1)
\,\,,
\]
and we get the inductive relation
$|\Re\, z_{k+1}|-1\geq  \exp^{\circ k}(|\Re\, z_1|-1)$.

Now if $\eta$ is sufficiently large, then indeed there exist $T$ and $t$ with
$T>t>0$ such that
\[
|\Re\, z_{k+1}|>\exp^{\circ k}(T)>
\exp^{\circ k}(t)> |\Im\, z_{k+1}|
\]
for almost all $k$. Once $k$ is so large that
$\exp^{\circ k}(T)>p\exp^{\circ k}(t)$, we have
$\exp^{\circ (k+1)}(T)>(\exp^{\circ (k+1)}(t))^p$. 
The assertion of the lemma follows.
\qed

We can finally prove that the set $R$ of dynamic rays has
Hausdorff dimension $1$:

\begin{corollary}[Hausdorff Dimension of the Union of Dynamic Rays]
\label{CorHausdorffRays} \lineclear
The set $R$ consisting of all dynamic rays has Hausdorff
dimension $1$.
\end{corollary}
\proof
Consider an arbitrary point $z$ of $R$, say $z=g(t)$ for some ray $g$ and some
$t>0$. Let $w:=g(t')$ for some $t'$ in $(0,t)$. Then by
Theorem~\ref{ThmEscapingSine} there is an $h$ not exceeding $\pi$ such that
$|\Im(f^{\circ k}(z)-f^{\circ k}(w))|\leq h$ for all $k$, and
$|\Re f^{\circ k}(z)|-|\Re f^{\circ k}(w) |\to\infty$ as $k\to\infty$.%
\footnote{Strictly speaking, we have stated Theorem~\ref{ThmEscapingSine} 
only for certain maps $z\mapsto ae^z+be^{-z}$ as specified in the
theorem, and only such maps will be used in the following sections,
so one can read this entire paper with only the maps $z\mapsto k\sinh z$ in mind.
However, the results in this section are true for all maps $z\mapsto ae^z+be^{-z}$ 
with $a$ and $b$ in $\C\sm\{0\}$.}
Fix $p$ with $p>1$. For each choice of $\xi>0$  Lemma~\ref{LemHorizExpansion} implies that 
there is an $N$ such that $f^{\circ N}(z)$ lies in $I_{p,\xi}$.

We have thus shown that $R\subset\bigcup_{N\geq 0}
f^{-N}(I_{p,\xi})$. If $\xi$ is sufficiently large, Lemma~\ref{LemDimParabola} ensures that 
$\dim_H(I_{p,\xi})\leq 1+1/p$. Now for each $N$ the set
$f^{-N}(I_{p,\xi})$ is a countable union of holomorphic preimages
of $I_{p,\xi}$, so parts 2 and 5 of Theorem~\ref{ThmElementaryProperties} imply that
$\dim_H(f^{-N}(I_{p,\xi}))\leq 1+1/p$.  It follows that $\dim_H(R)\leq 1+1/p$.
Since this is true for every $p$ greater than $1$, we conclude that $\dim_H(R)\leq 1$. 
Equality follows because $R$ contains curves.
\qed

Now we have our dimension paradox complete for $f(z)=k\pi\sinh z$, using
Theorem~\ref{ThmEscapingSine} (which still requires proof): every point $z$ of $\C$ either lies
on a dynamic ray, and thus is in $R$, or it is a landing point of one
or several dynamic rays in $R$ that connect $z$ to $\infty$.
Since the set $R$ has Hausdorff dimension $1$ (hence planar Lebesgue measure
zero), the set $E=\C\sm R$ has full measure and is in fact
everything but the one-dimensional set $R$. This proves
Theorem~\ref{ThmParadox} (further details can be found in \cite{Cosinus}).

\begin{figure}[htbp]
\includegraphics[viewport=050 020 580 500,height=80mm,clip]{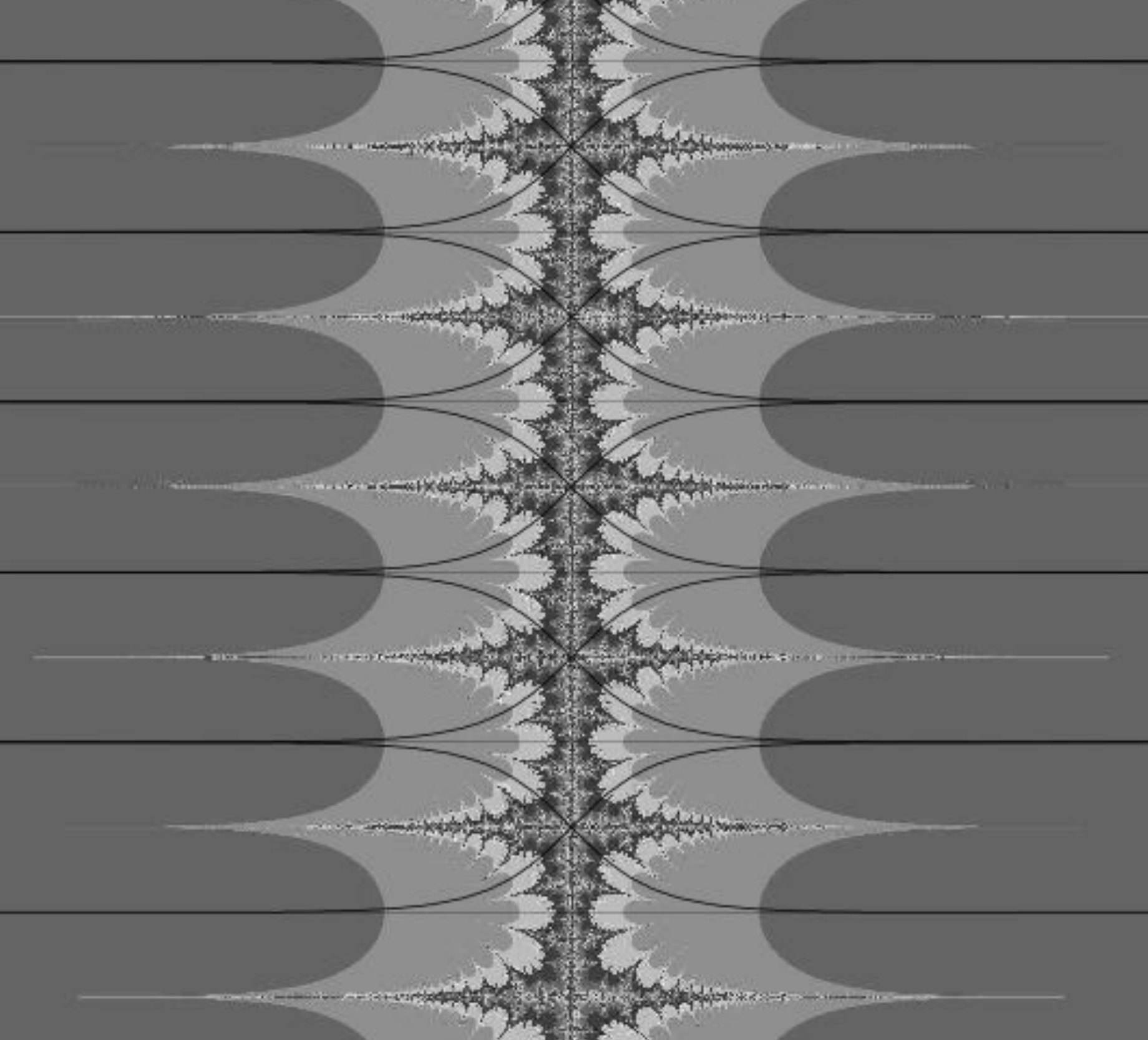}
\caption{The dynamical plane of the map $f\colon z\mapsto
\pi\sinh z$. Several dynamic rays are shown.
}
\label{FigSinDynamics}
\end{figure}

\MAAsection{DYNAMICAL FINE-STRUCTURE OF THE HYPERBOLIC SINE MAP}
\label{SecFineStructure}
We now proceed to explain why Theorem~\ref{ThmEscapingSine}
is true, and why it is interesting from the perspective of
dynamical systems. For simplicity, we restrict attention to maps
$f(z)=k\pi\sinh z=(k\pi/2)(e^z-e^{-z})$ with $k$ a positive integer (see
Figure~\ref{FigSinDynamics}).

First observe that $f$ is periodic with period $2\pi i$ ($f$ is the rotated
sine function) and maps $i\R$ onto the interval
$[-k\pi i,k\pi i]$. Notice also that $f\colon\R\to\R$
is a homeomorphism with $f(0)=0$ and $f'(x)\geq\pi$ for all
$x$ in $\R$, from which it follows that $\R\sm\{0\}$ is contained in 
the escape set $I$. In fact, $\R^+$
and $\R^-$ are two of the path components of $I$: they are both
dynamic rays, and they connect each of their points to $\infty$
through $I$. Since $f(z+i\pi)=-f(z)$, other dynamic rays include the
curves $i\pi n+\R^+$ and $i\pi n+\R^-$ for integers $n$. These map
under $f$ onto $\R^+$ or $\R^-$. This gives a useful partition
for the dynamics: for $n$ in $\Z$ set 
\begin{eqnarray*}
U_{n,R}
&:=&
\{z\in\C\colon \Re\, z >0, \Im\,z\in(2\pi n, 2\pi (n+1))\}
\,\,,
\\
U_{n,L}
&:=&
\{z\in\C\colon \Re\,z <0, \Im\,z\in(2\pi n, 2\pi (n+1))\}
\,\,.
\end{eqnarray*}
(This is an ad hoc partition for our special maps $f$ that uses
the symmetry given by the invariant real and imaginary axes. In
\cite{Cosinus}, a different partition is used that works for
more general maps $f$.)

The geometry of the mapping $f$ is such that its restrictions are conformal isomorphisms 
\[
f\colon U_{n,R}\to \C\sm(\R^+\cup[-k\pi i,k\pi i])
\] 
and 
\[
f\colon U_{n,L}\to \C\sm(\R^-\cup[-k\pi i,k\pi i])
\,\,,
\]
so the image of each
$U_{n,\times}$ is a one-sheeted covering of $U_{n,\times}$.
This is a useful property, called the {\em Markov property}, that aids in reducing
many dynamical questions to questions about symbolic dynamics. 

Let $\Z_R:=\{\dots,-2_R,-1_R,0_R,1_R,
2_R,\dots\}$ and $\Z_L:=\{\dots,-2_L$, $ -1_L$, $0_L,1_L,
2_L,\dots\}$ be two disjoint copies of $\Z$, and let
$\Sym\colon= (\Z_R\cup\Z_L)^\N$ be the space of sequences with elements in
$\Z_R\cup\Z_L$. To each $z$ in $\C$ we assign an {\em itinerary}
$\s=s_1s_2s_3\dots$ in $\Sym$ such that $s_k=n_R$ if $f^{\circ
(k-1)}(z)$ is in $\ovl U_{n,R}$ and $s_k=n_L$ if $f^{\circ
(k-1)}(z)$ is in $\ovl U_{n,L}$. There are ambiguities if the orbit of
$z$ ever enters $\R$ or $[-k\pi i,k\pi i]$, but such points are
easy to understand anyway, and we admit all itineraries in such
cases (the number of possible itineraries for a given point $z$ can be 
as large as four; see the discussion in the proof of Theorem~\ref{ThmEscapingSine}).
The following lemma furnishes a mechanism for understanding the
detailed dynamics of $f$:

\begin{lemma}[Symbolic Dynamics and Curves]
\label{LemCurves} \lineclear
For each sequence $\s$ in $\Sym$ the set of all points $z$ in $\C$
with itinerary $\s$ is either empty or a curve that connects
$\infty$ to a well-defined landing point in $\C$. 
For each such curve each of its points other than the landing point escapes.
\end{lemma}
\sketch
For each positive $N$ let $U_{\s,N}$ be the set of points $z$ such that
the first $N$ entries in the itinerary of $z$ coincide with the first $N$ entries
of $\s$. With the aid of the Markov property it is quite easy to see that
each $\ovl U_{\s,N}$ is
a closed, connected, and unbounded subset of $\C$. Moreover,  in the topology
of the Riemann sphere, adding the point $\infty$ to these sets
yields compact and connected sets containing $\infty$. Let
\[
C_{\s}:=\bigcap_{N\in\N} (\ovl U_{\s,N}\cup\{\infty\})
\,\,.
\]
This is obviously a nested intersection, so $C_{\s}$ is compact and
connected and contains $\infty$. If $C_{\s}=\{\infty\}$,
then we have nothing to prove. Otherwise, we can show that $f$ is
expanding enough so that for any two points $z$ and $w$ in $C_{\s}$ and any
$\eta>0$ there is an $n$ such that $\left||\Re f^{\circ
n}(z)|-|\Re f^{\circ n}(w)|\right|>\eta$.
Lemma~\ref{LemHorizExpansion} implies then that at least one of the points 
$z$ and $w$ escapes. (The
expansion comes from the fact that $U:=\C\sm \{-i\pi,0,i\pi\}$
carries a unique normalized hyperbolic metric and that
$f^{-1}(U)\subset U$. With respect to this metric on $U$, every
local branch of $f^{-1}$ is contracting, which makes $f$ locally
expanding. This argument requires nothing but the fact that the
universal cover of $U$ is $\disk$, plus the Schwarz lemma on
holomorphic self-maps of $\disk$.)

It follows that all points in $C_\s\sm\{\infty\}$ escape, with at most
one exception; the estimates in Lemma~\ref{LemHorizExpansion} imply
that these points escape extremely fast. This means that for almost
all $z$ in $C_\s\sm\{\infty\}$ we have $f^{\circ
n}(z)\to\infty$ very fast, hence $|(f^{\circ n})'(z)|\to\infty$
very fast. Thus the forward iterates of $z$ are very strongly
expanding. Conversely, if $z_n:=f^{\circ n}(z)$, then the branch
of $f^{-n}$ sending $z_n$ to $z$ is strongly contracting. This
implies that the boundaries of the $U_{\s,N}$, which are curves,
converge locally uniformly to $C_\s$. This ensures that $C_\s$
is a curve.
\qed

This lemma is all we need to establish the two main results about the
dynamics of the function $f$.

\proofof{Theorem~\ref{ThmEscapingSine}}
Every point $z$ in $\C$ has at least one associated itinerary. If it
has more than one, then under iteration it must map into $i\R$ or
into $\R+2\pi i\Z$. In the latter case, the next iteration lands
in $\R$, so the orbit reaches either the fixed point $0$ or one
of the two dynamic rays $\R^+$ or $\R^-$. If the orbit reaches
$i\R$, then from that iteration on it spends its entire forward orbit in the
interval $[-k\pi i,k\pi i]$; in particular, the orbit is
bounded. Therefore, a point has four itineraries if and only if
its orbit eventually terminates at $0$. A point has two
itineraries if it lands in the invariant interval $[-k\pi i,k\pi i]$
(and has bounded orbit), or if it lands in $\R^+\cup\R^-$ and
escapes. Every  other point has a single itinerary.

Recall that the set of points with a given itinerary is a single
dynamic ray consisting of escaping points, together with the unique
landing point of the ray (Lemma~\ref{LemCurves}). This implies
that every point in $\C$ either lies on a unique dynamic ray or is the
landing point of one, two, or four dynamic rays. This proves
statements~\ref{ItemThm2} and \ref{ItemThm3} in the theorem.

For statement~\ref{ItemThm1}, we have constructed rays
consisting of escaping points, and the partition makes it clear
that every ray has real parts tending to $\pm\infty$, while the
imaginary parts are constrained to some interval of length $\pi$. It
is clear that each escaping point either is on a unique ray or is
the landing point of a ray; if a ray lands at an escaping
point, then the landing point neither lies on any other ray nor is the
landing point of another ray. Therefore, each ray (possibly
together with its endpoint) is contained in a path component of
$I$. It is also true that each path component of $I$ consists of
a single ray, possibly together with its endpoint. The proof of this fact
requires some ingredients from continuum theory
(see \cite[sec.~4]{FRS}).
\qed

\MAAsection{LEBESGUE MEASURE AND ESCAPING POINTS}
From the point of view of dynamical systems, an important
question to ask is the following: {\em What do most orbits do under iteration?}
From a topological vantage point, most points are on dynamic
rays, rather than being endpoints of rays. On the other hand, since
the union of the rays has Hausdorff dimension $1$, measure theory
says that most points are endpoints of rays. However, as we will now see, even
measure theory asserts that most points in $\C$ escape (for our maps $z\mapsto k\pi\sinh z$):
this assertion is a combination of results of McMullen~\cite{McMullen} and
Bock~\cite{Bock}. As a result, almost all
points are escaping endpoints of rays. Along the way, we visit a result of
Schubert~\cite{Hendrik} that settles a conjecture of
Milnor~\cite[sec.~6]{MiIntro} in the affirmative.

%\pagebreak

\begin{theorem}[Lebesgue Measure of Escaping Points]
\label{ThmMostPointsEscape} \lineclear
\begin{enumerate}
\item
For every map $z\mapsto\lambda e^z$ with $\lambda\neq 0$ the
set $I$ of escaping points has two-dimensional Lebesgue measure zero but
Hausdorff dimension $2$ \cite{McMullen}.
\item
However, for every map $z\mapsto ae^z+be^{-z}$ with $ab\neq 0$
the set $I$ has infinite two-dimensional Lebesgue
measure \cite{McMullen}.
For every strip $S=\{z\in\C\colon \alpha\le\Im\,z\le\beta\}$ in $\C$
the two-dimensional Lebesgue measure of $S\sm I$ is finite
\cite{Hendrik}.
\end{enumerate}
\end{theorem}
\sketch
Choose $\xi>0$, and set
$\half_\xi=\{z\in\C\colon\Re\,z>\xi\}$. We show that for every
map $E(z)=\lambda\exp(z)$ and sufficiently large $\xi$ the set
$Z_\xi:=\{z\in\C\colon\Re\, E^{\circ n}(z)>\xi \mbox{ for all
$n$}\}$ has measure zero. In fact, for each square $Q$ in $\half_\xi$ of side-length
$2\pi$ with sides parallel to the coordinate axes the image $E(Q)$ is a
large annulus in $\C$, and the probability that a point $z$ in $Q$ has
$E(z)$ in $\half_\xi$ is approximately $1/2$. The chance of
surviving $n$ consecutive iterations in $Z_\xi$ is then $2^{-n}$
(assuming independence of probabilities in the consecutive steps).
Hence, the set of points $z$ in $Q$ whose entire orbits lie in
$\half_\xi$ has two-dimensional Lebesgue measure zero, and thus 
all of $Z_\xi$ has two-dimensional Lebesgue measure
zero. But since $|E(z)|=|\lambda|\exp(\Re\,z)$, for every point $z$ in
$I$ there must be an $N$ such that $E^{\circ n}(z)$ belongs to  $Z_\xi$ for {\em
all} $n$ with $n\geq N$. Since $I$ is a subset of $\bigcup_{n\geq 0}E^{\circ -n}(Z_\xi)$,
it has measure zero for exponential maps $z\mapsto \lambda e^z$.

The situation is different for $E(z)=ae^z+be^{-z}$ with $ab\neq 0$: instead of
throwing away half of the points in every step, we can
``recycle'' (in a literal sense) most of them: this time
$|E^{\circ n}(z)|\to\infty$ implies that $|\Re\, E^{\circ
n}(z)|\to\infty$. We use another parabola (or rather the
complement thereof), namely,
\[
P:=\{x+iy\in\C\colon |y|<|x|^2\} \,\,.
\]
If $z=x+iy$ with $|x|$ sufficiently large is such that $E(z)$ lies in
$P$, then 
\[
|\Re\,E(z)| \geq |E(z)|^{1/2} \approx e^{|x|/2}\gg |x|
\,\,,
\]
so points that escape to $\infty$ within $P$ do so quite rapidly. On the other hand, the
image of a square $Q$ as in the first part (with real parts $x$ greater than
$\xi$ or less than $-\xi$) is again an annulus, but the fraction
of $E(Q)$ within $P$ is approximately
$1-e^{-|x|/2}$, so most of the points survive the first step.
Among these, a fraction of $1-e^{-(e^{-|x|/2})/2}$ survives the
second step, and so on.  The total fraction of points within $Q$
that ``get lost'' from
$P$ under iteration is less than $1$, from which we infer that $I\cap Q$ has
positive two-dimensional Lebesgue measure \cite{McMullen}. To be more
precise, we recursively define a sequence $(\xi_n)$ by $\xi_0=\xi$ and
$\xi_{n+1}=e^{|\xi_n|/2}$ for $n=1,2,\dots$. Then  $\xi_n>2^{n-1}\xi_1$
for all $n$, provided that $\xi_0$ is sufficiently large.
If
\[
Q_n:=\{z\in
Q\colon E^{\circ k}(z)\in P\mbox{ for $k=0,1,2,\dots,n$}\}
\,\,,
\]
then each $z$ in $Q_n$ has $|\Re\, E^{\circ n}(z)|>\xi_n$. This means that of all 
the points in $Q_n$, a fraction of at least 
$1-e^{-\xi_n/2}=1-1/\xi_{n+1}$ survives one more iteration within $P$.
Thus, denoting two-dimensional Lebesgue measure by $\mu$, we get
\[
\frac{\mu(Q_n)}{\mu(Q)}
> 1-\frac{1}{\xi_1}-\frac{1}{\xi_2}-\dots\frac{1}{\xi_n}
> 1-\frac{2}{\xi_1}
=1-2e^{-\xi_0}
\,\,.
\]
Since $\bigcap_n Q_n$ is contained in $I$, it follows that 
\[
\mu(I\cap Q)>(1-2e^{-\xi_0/2})\mu(Q) 
\,\,:
\]
the set of escaping points has positive density in $Q$, hence
$I$ has positive (even infinite) two-dimensional Lebesgue measure.

In fact, we have shown much more:
${\mu(Q\sm I)}<2e^{-\xi_0/2}{\mu(Q)}$ \cite{Hendrik}. Therefore,
for each horizontal strip $S$ of height $2\pi$ the complement of
$I$ in $S$ has finite Lebesgue measure: 

\[
\mu\left(z\in S\sm I\colon |\Re\,z|>\xi_0\right)
<2\pi \int_{\xi_0}^\infty 2e^{-x/2}\,dx
=
4\pi e^{-\xi_0/2}
\,\,.
\]
This proves the result.
\qed

\remark
Milnor~\cite[sec.~6]{MiIntro}  conjectured that for
$f(z)=\sin z$ the set of points converging to the fixed point
$z=0$ has finite Lebesgue area in every strip $S'=\{z\in\C\colon \alpha\le\Re\,z\le\beta\}$. 
Since sine and hyperbolic sine represent the same map in rotated
coordinate systems, this follows from Schubert's result.

We conclude with another special case in which the set $I$ is so
large that $\C\sm I$ has measure zero \cite{Cosinus}:

\begin{corollary}[Escaping Set of Full Measure]
\lineclear
For maps $z\mapsto k\pi\sin z$ or $z\mapsto k\pi\sinh z$ with a nonzero integer $k$
the set $I\cap E$ has full two-dimensional Lebesgue measure 
{\em (}i.e., the measure of $\C\sm (I\cap E)$ is zero{\em).}
\end{corollary}
\proof
We invoke a theorem of Bock~\cite{Bock}: for
an arbitrary transcendental entire function at least one of the
following two statements holds: (i) almost every orbit is dense
in $\C$ or (ii) almost every orbit converges to $\infty$ or to
one of the critical orbits (a critical orbit is the orbit
of one of the two critical values $\pm k\pi i$). But since $I$ has
positive measure, case (i) cannot hold, so statement (ii)
follows. For the map $E:z\mapsto k\pi\sinh z$ (or, equivalently, 
$z\mapsto k\pi\sin z$) the two critical values map to the fixed point
$0$. However, since $|E'(0)|>1$ (i.e., the fixed point $0$ is
``repelling''), the only points whose orbits can converge to $0$ are
those countably many points that land exactly on $0$ after finitely
many iterations. Therefore, almost every orbit must escape.
\qed

\noindent
{\bf ACKNOWLEDGMENTS.} I would like to express my gratitude to
Bogusia Karpi\'nska for many interesting discussions and for
allowing me to include her example in section~\ref{SecKarpinska}.
I would also like to thank Cristian Leordeanu and G\"unter 
Rottenfu{\ss}er for their help with the illustrations in this article.

\medskip

\hide{
\noindent
{\bf DIERK SCHLEICHER} studied physics and mathematics in Hamburg and Cornell, with semesters abroad in Princeton and Paris. He held teaching and research positions in Munich and Stony Brook, before moving to Bremen in 2001 to help build up International University Bremen: a small and exciting new university with strong students from 80 countries. He has spent research semesters in Berkeley, Paris, and Toronto. His research interests are on the interplay between geometry and dynamics, with a focus on holomorphic dynamics and symbolic dynamics. On campus, he enjoys spending time with his student advisees and organizing events like the International Mathematics Olympiad 2009 in Bremen; off campus, he likes kayaking, paragliding, and hiking in the mountains. 
}

\end{document}